\documentclass[12pt]{article}
\pdfoutput=1 
\usepackage[T1]{fontenc}
\usepackage[active]{srcltx}
\usepackage{lmodern}
\usepackage{amsmath,amstext,amssymb,graphicx,graphics,color}
\usepackage{subcaption}
\usepackage{tikz}
\usepackage{pgfplots}
\usepackage{theorem}
\usepackage{cite}               
\usepackage{hyperref}           
\usepackage{bbm}                
\usepackage{fancybox}           
\usepackage[section]{placeins}  
\theorembodyfont{\normalfont\slshape} 
\newtheorem{definition}{Definition}
\newtheorem{theorem}{Theorem}
\newtheorem{proposition}{Proposition}[section]
\newtheorem{corollary}[proposition]{Corollary}

\newtheorem{lemma}[proposition]{Lemma}
\theoremstyle{break} 

{{\par\noindent \bf Proof. \nobreak}}%
{\nobreak \removelastskip \nobreak \hfill $\Box$ \medbreak}
{{\par\noindent \bf Proof \nobreak}}%
{\nobreak \removelastskip \nobreak \hfill $\Box$ \medbreak}
{{\par\noindent \bf Proof lemma. \nobreak}}%
{\nobreak \removelastskip \nobreak \bf End proof lemma. \medbreak}
\newenvironment{remark}{\par \medskip \noindent {\bf Remark. }\nobreak}{\par \medskip}
\usepackage{fancyhdr}
 \fancypagestyle{plain}{
  \fancyhead{}
  
}

\def\paragraph#1{{\bf #1\ }}
\newcommand{\expo}{\mathrm{e}}

\newcommand{\dd}{\mathrm{d}}
\newcommand{\EE}{\mathrm{E}}
\newcommand{\R}{{\mathbb{R}}}
\newcommand{\eps}{\varepsilon}
\newcommand{\emp}{\mathrm{emp}}
\newcommand{\NN}{\mathrm{N}}
\newcommand{\HH}{\mathrm{H}}
\newcommand{\overbar}[1]{\mkern 1.5mu\overline{\mkern-1.5mu#1\mkern-1.5mu}\mkern 1.5mu}
\usepackage[utf8]{inputenc}
\usepackage{unicode_character}

\def\Proof{\noindent{\bf Proof}\quad}
\def\qed{\hfill$\square$\smallskip}

\title{Entropy dissipation and propagation of chaos for the uniform reshuffling model}
\author{Fei Cao\footnotemark[1] \and Pierre-Emmanuel Jabin\footnotemark[2] \and Sebastien Motsch\footnotemark[1]}

\begin{document}

\maketitle

\footnotetext[1]{Arizona State University - School of Mathematical and Statistical Sciences, 900 S Palm Walk, Tempe, AZ 85287-1804, USA}
\footnotetext[2]{Pennsylvania State University - Department of Mathematics and Huck Institutes, State College, PA 16802, USA}

\tableofcontents

\begin{abstract}
We investigate the uniform reshuffling model for money exchanges: two agents picked uniformly at random redistribute their dollars between them. This stochastic dynamics is of mean-field type and eventually leads to a exponential distribution of wealth. To better understand this dynamics, we investigate its limit as the number of agents goes to infinity. We prove rigorously the so-called propagation of chaos which links the stochastic dynamics to a (limiting) nonlinear partial differential equation (PDE). This deterministic description, which is well-known in the literature, has a flavor of the classical Boltzmann equation arising from statistical mechanics of dilute gases. We prove its convergence toward its exponential equilibrium distribution in the sense of relative entropy.
\end{abstract}



\section{Introduction}


Econophysics is an emerging branch of statistical physics that incorporate notions and techniques of traditional physics to economics and finance \cite{savoiu_econophysics:_2013,chatterjee_econophysics_2007,dragulescu_statistical_2000}. It has attracted considerable attention in recent years raising challenges on how various economical phenomena could be explained by universal laws in statistical physics, and we refer to \cite{chakraborti_econophysics_2011,chakraborti_econophysics_2011-1,pereira_econophysics_2017,kutner_econophysics_2019}  for a general review.

The primary motivation for study models arising from econophysics is at least two-fold: From the perspective of a policy maker, it is important to deal with the raise of income inequality \cite{dabla-norris_causes_2015,de_haan_finance_2017} in order to establish a more egalitarian society. From a mathematical point of view, we have to understand the fundamental mechanisms, such as money exchange resulting from individuals, which are usually agent-based models. Given an agent-based model, one is expected to identify the limit dynamics as the number of individuals tends to infinity and then its corresponding equilibrium when run the model for a sufficiently long time (if there is one), and this guiding approach is carried out in numerous works across different fields among literature of applied mathematics, see for instance \cite{naldi_mathematical_2010,barbaro_phase_2014,carlen_kinetic_2013}.

\begin{figure}[!h]
  \centering
  \includegraphics[scale=1]{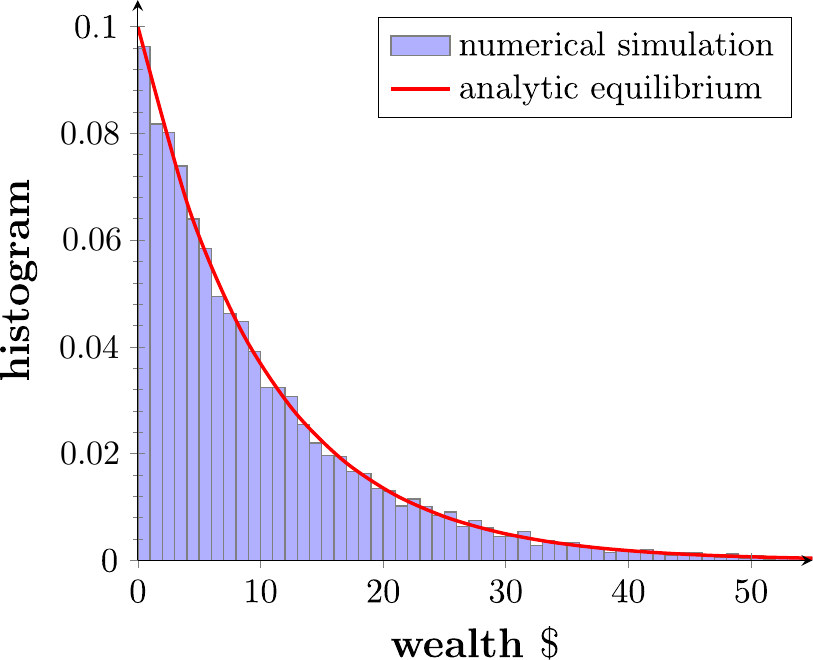}
  \caption{Simulation results for the uniform reshuffling model. The blue histogram shows the distribution of money after $T=1000$ time unit. The red solid curve is the limiting exponential distribution proved in \cite{lanchier_rigorous_2018}. We used $N=10,000$ agents, each starting with $\$10$.}
  \label{Simul} 
\end{figure}

In this work, we consider the so-called uniform reshuffling model for money exchange in a closed economic system with $N$ agents. The dynamics consists in choosing at random time two individuals and to redistribute their money between them. To write this dynamics mathematically, we denote by $X_i(t)$ the amount of dollar agent $i$ has at time $t$ for $1\leq i\leq N$. At a random time generated by a Poisson clock with rate $N$, two agents (say $i$ and $j$) update their purse according to the following rule:
\begin{equation}
\label{uniform}
\big(X_i,X_j\big) \leadsto \big(U(X_i\!+\!X_j)\,,\,(1\!-\!U)(X_i\!+\!X_j)\big),
\end{equation}
where $U$ is a uniform random variable over the interval $[0,1]$ (i.e. $U \sim \mathrm{Uniform}[0,1]$). The \emph{uniform reshuffling model} is first studied in \cite{dragulescu_statistical_2000} via simulation. The agent-based numerical simulation suggests that, as the number of agents and time go to infinity, the limiting distribution of money approaches the exponential distribution as shown in Figure \ref{Simul}.

It is well-known (see for instance \cite{matthes_steady_2008,bassetti_explicit_2010,during_kinetic_2008,apenko_monotonic_2013}) that under the large population $N \to \infty$ limit, We can formally show that the law of the wealth of a typical agent (say $X_1$) satisfies the following limit PDE in a weak sense:
\begin{equation}\label{uniform_reshuffling_pde}
\partial_t q(t,x) = \int_0^\infty\int_0^\infty \frac{\mathbbm{1}_{[0,k+\ell]}(x)}{k+\ell}q(t,k)q(t,\ell)\,\dd \ell \, \dd k - q(t,x).
\end{equation}
To our best knowledge, the rigorous derivation of the limit equation \eqref{uniform_reshuffling_pde} from the particle system description is absent in most of the literature on econophysics (just like many other PDEs arising from models in econophysics \cite{katriel_immediate_2015,heinsalu_kinetic_2014,chakrabarti_econophysics_2013}), because the propagation of chaos effect is implicitly assumed in the large $N$ limit in most derivations. The remarkable exception is the paper \cite{cortez_particle_2018}, where the author showed a uniform-in-time propagation of chaos by virtue of a delicate coupling argument. In section \ref{sec:propagation of chaos} of this manuscript, we will provide an alternative rigorous justification of the equation \eqref{PDE} under the limit $N \to \infty$.

\begin{figure}[!htb]
\centering
\includegraphics[scale=0.9]{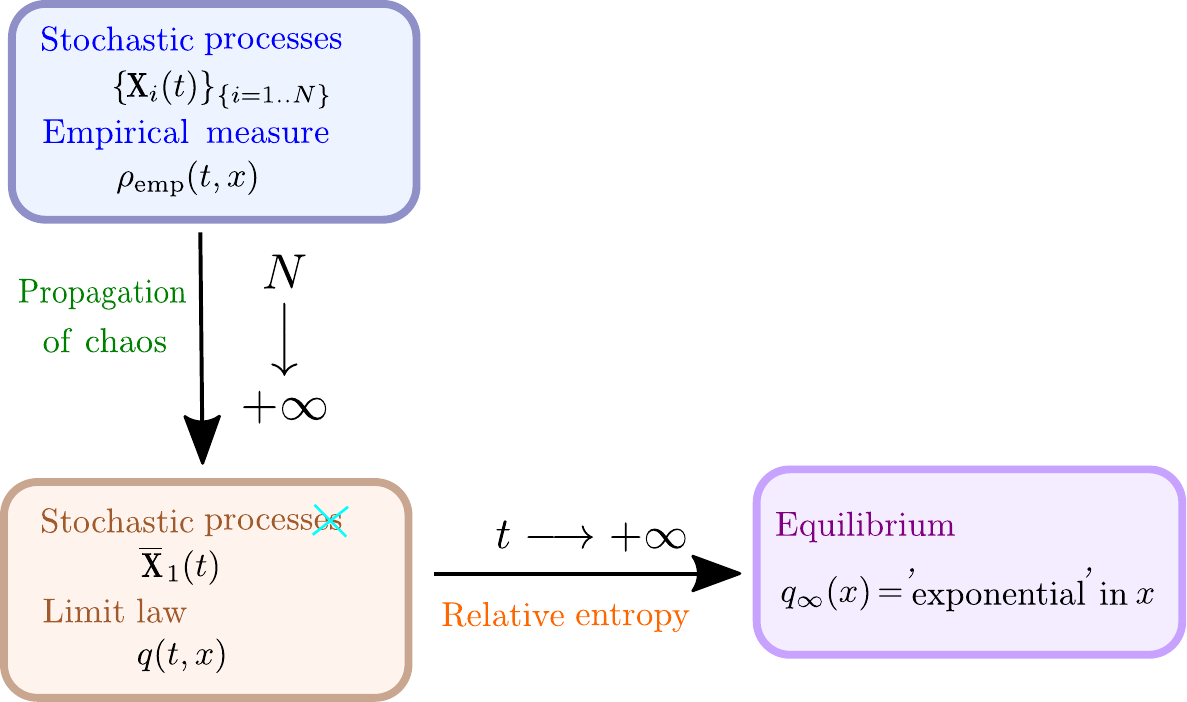}
\caption{Schematic illustration of the general strategy of our treatment of the uniform reshuffling dynamics.}
\label{schema} 
\end{figure}

Once the limit PDE is identified from the interacting particle system, the natural next step is to study the problem of convergence to equilibrium of the PDE at hand, it has been shown in \cite{during_kinetic_2008,matthes_steady_2008} that the unique (smooth) solution of \ref{uniform_reshuffling_pde} converges to its exponential equilibrium distribution exponentially fast in Wasserstein and Fourier metrics. In the present work, we demonstrate a polynomial in time convergence in relative entropy, by establishing a entropy-entropy dissipation inequality (see Theorem \ref{EEP} below) which is not available among the literature. An illustration of the general strategy used in this work (and implicitly in many of the works cited above) is shown in Figure \ref{schema}.

Although only a very specific binary exchange model is explored in the present paper, other exchange rules can also be imposed and studied, leading to different models. To name a few, the so-called immediate exchange model introduced in \cite{heinsalu_kinetic_2014} assumes that pairs of agents are randomly and uniformly picked at each random time, and each of the agents transfer a random fraction of its money to the other agents, where these fractions are independent and uniformly distributed in $[0,1]$. The uniform reshuffling model with saving propensity investigated in \cite{chakraborti_statistical_2000,lanchier_rigorous_2018} suggests that the two interacting agents keep a fixed fraction $\lambda$ of their fortune and only the combined remaining fortune is uniformly reshuffled between the two agents, which makes the uniform reshuffling model the particular case $\lambda = 0$. For more variants of exchange models with (random) saving propensity and with debts, we refer the readers to \cite{chatterjee_pareto_2004} and \cite{lanchier_rigorous_2018-1}.

This manuscript is organized as follows: in section \ref{sec:limit_pde}, we briefly discuss the properties of the limit equation \eqref{uniform_reshuffling_pde}. We show in section \ref{sec:wasserstein_linearization} convergence results for the solution of \eqref{uniform_reshuffling_pde} in Wasserstein distance and in the linearized region. We take on the most delicate analysis of the entropy-entropy dissipation relation in section \ref{sec:entropy_dissipation}. Finally, we present a rigorous treatment of the propagation of chaos phenomenon in section \ref{sec:propagation of chaos}.

\section{The limit PDE and its properties}
\label{sec:limit_pde}
\setcounter{equation}{0}

We present a heuristic argument behind the derivation of the limit PDE \eqref{PDE} arising from the uniform reshuffling dynamics in section \ref{subsec:2.1}. Several elementary properties of the solution of \eqref{PDE} are recorded in section \ref{subsec:2.2}. Section \ref{subsec:2.3} is devoted to another formulation of the uniform reshuffling model, which can be viewed as a \emph{lifting} of the reshuffling mechanics \eqref{uniform} and is implicitly exploited in \cite{apenko_monotonic_2013}. In section \ref{subsec:2.4}, we highlight a key ingredient known as the \emph{micro-reversibility}, of the collision operator determined by the right side of \eqref{PDE}, which allows us to construct certain Lyapunov functions associated with \eqref{PDE} (such as entropy).

\subsection{Formal derivation of the limit PDE}
\label{subsec:2.1}

Introducing $\mathrm{N}^{(i,j)}_t$ independent Poisson processes with intensity $1/N$, the dynamics can be written as:
\begin{equation}
  \label{eq:SDE_uniform}
  \dd X_i(t) = \sum_{j=1..N,j\neq i} \Big(U(t-)(X_i(t-)\!+\!X_j(t-)) - X_i(t-)\Big) \dd \mathrm{N}^{(i,j)}_t
\end{equation}
with $U(t)\sim \mathrm{Uniform}[0,1]$ independent of $\{X_i(t)\}_{1\leq i\leq N}$. As the number of players $N$ goes to infinity, one could expect that the processes $X_i(t)$ become independent and of same law. Therefore, the limit dynamics would be of the form:
\begin{equation}
  \label{eq:SDE_uniform_limit}
  \dd \overline{X}(t) = \Big(U(t-)(\overline{X}(t-)\!+\!\overline{Y}(t-)) - \overline{X}(t-)\Big) \dd \overline{\mathrm{N}}_t
\end{equation}
where $\overline{Y}(t)$ is an independent copy of $\overline{X}(t)$ and $\overline{\mathrm{N}}_t$ a Poisson process with intensity $1$. Taking a test function $φ$, the weak formulation of the dynamics is given by:
\begin{equation}
  \label{eq:SDE_uniform_limit_weak}
  \dd 𝔼[φ(\overline{X}(t))] = 𝔼[φ\Big(U(t)(\overline{X}(t)\!+\!\overline{Y}(t))\Big) - φ(\overline{X}(t))]\,\dd t
\end{equation}
In short, the limit dynamics correspond to the jump process:
\begin{equation}
  \label{uniform_limit}
  \overline{X} \leadsto U(\overline{X}\!+\!\overline{Y}).
\end{equation}
Let's denote $q(t,x)$ the law of the process $\overline{X}(t)$. To derive the evolution equation for $q(t,x)$, we need to translate the effect of the jump of $\overline{X}(t)$ via \eqref{uniform_limit} onto $q(t,x)$.

\begin{lemma}(Hierarchy of probability distributions)\label{lem1}
  Suppose $X$ and $Y$ two independent random variables with probability density $q(x)$ supported on $[0,\infty)$. Let $Z = U(X+Y)$ with $U\sim \mathrm{Uniform}([0,1])$ independent of $X$ and $Y$. Then the law of $Z$ is given by $Q_+[q]$ with:
  \begin{eqnarray}
    \label{eq:Q_plus}
    Q_+[q](x) &=& \int_{m=0}^\infty \frac{\mathbbm{1}_{[0,m]}(x)}{m}\left(\int_{z=0}^m q(z)q(m-z)\dd z\right)\dd m\\
              &=& \int_{\mathbb{R}_+ \times \mathbb{R}_+} \frac{\mathbbm{1}_{[0,k+\ell]}(x)}{k+\ell}q(k)q(\ell)\dd \ell \dd k
  \end{eqnarray}
\end{lemma}
\Proof Let's introduce a test function $φ$.
\begin{eqnarray*}
  𝔼[φ(U(X\!+\!Y))] &=& ∫_{x \geq 0} ∫_{y \geq 0} ∫_{u=0}^1 φ(u(x+y))q(x)q(y)\,\dd u\dd x\dd y \\
                   &=& ∫_{m \geq 0} ∫_{z = 0}^m ∫_{u=0}^1 φ(um)q(z)q(m-z)\,\dd u\dd z\dd m \\
                   &=& ∫_{m \geq 0} ∫_{z = 0}^m \frac{1}{m}∫_{s=0}^m φ(s)q(z)q(m-z)\,\dd s\dd z\dd m
\end{eqnarray*}
using the change of variables $z=x$ and $m=x+y$ followed by $s=um$. We conclude using Fubini that:
\begin{eqnarray}
  \label{eq:weak_Q_plus}
  𝔼[φ(U(X\!+\!Y))] &=& ∫_{s\geq 0} φ(s)  \left(∫_{m \geq 0} \mathbbm{1}_{[0,m]}(s)  \frac{1}{m} ∫_{z = 0}^mq(z)q(m-z)\,\dd z\dd m \right)\dd s \nonumber \\
  &=& ∫_{s\geq 0} φ(s) Q_+[q](s)\,\dd s
\end{eqnarray}
with $Q_+[q]$ defined by \eqref{eq:Q_plus}.
\qed



We can now write the evolution equation for the law of $\overline{X}(t)$ \eqref{eq:SDE_uniform_limit}, the density $q(t,x)$ satisfies weakly:
\begin{equation}
  \label{PDE}
  \partial_t q(t,x) = G[q](t,x) \qquad \text{ for } t\geq 0 \text{ and } x\geq 0
\end{equation}
with
\begin{equation}
  \label{G}
  G[q](x) := Q_+[q](x)-q(x)
  = \int_0^\infty\int_0^\infty \frac{\mathbbm{1}_{[0,k+\ell]}(x)}{k+\ell}q(k)q(\ell)\dd \ell \dd k - q(x).
\end{equation}

\subsection{Evolution of moments}
\label{subsec:2.2}
Now we will establish several elementary properties of the solution of \eqref{PDE}:

\begin{proposition}\label{conservation}
  Assume that $q(t,x)$ is a classical (and global in time) solution of \eqref{PDE} and define by $m_k(t)$ the $k$-th moment of $q$:
  \begin{equation}
    \label{eq:k_th_moment}
    m_k(t) :=  \int_0^\infty x^k q(t,x) \dd x.
  \end{equation}
  Then:
  \begin{equation}
    \label{eq:m_k}
    m_k'(t) = \frac{1}{k+1} ∑_{j=0}^k C_k^jm_j(t) m_{k-j}(t)  \;-\; m_k(t),
  \end{equation}
  where $C_k^j=\left(
      \begin{array}{c}
        k\\
        j
      \end{array}
\right)=\frac{k!}{j!(k-j)!}$ represents the binomial coefficient.
\end{proposition}
\Proof
Notice that the moment can be written as: $m_k(t) = 𝔼[\overline{X}^k(t)]$. Thus, we use the weak formulation of the evolution equation of $q(t,x)$ \eqref{eq:SDE_uniform_limit_weak} with $φ(x)=x^k$ and deduce that:
\begin{displaymath}
  m_k' = 𝔼[\big(U(\overline{X}\!+\!\overline{Y})\big)^k - \overline{X}^k]   \,= \, 𝔼[U^k]𝔼[(\overline{X}\!+\!\overline{Y}))^k] - m_k,
\end{displaymath}
since $U$ is independent of $\overline{X}$ and $\overline{Y}$. Moreover, $𝔼[U^k]=∫_{u=0}^1 u^k\,\dd u = \frac{1}{k+1}.$ Using the independence of $\overline{X}$ and $\overline{Y}$ and expanding lead to \eqref{eq:m_k}.
\qed

\begin{corollary}
  Let $q(t,x)$ solution of \eqref{PDE} and $m_k(t)$ its $k-$th moment \eqref{eq:k_th_moment}. The total mass and the mean are preserved, i.e. $m_0'(t)=m_1'(t)=0$ and all the moments $m_k(t)$ converges in time exponentially fast.
\end{corollary}
\Proof
Writing \eqref{eq:m_k} for $k=2$ leads to:
  \begin{equation}
    \label{eq:m_2}
    m_2' = -\frac{1}{3} m_2 + \frac{2}{3} m_1^2
  \end{equation}
  and thus $m_2(t) = 2\,m_1^2 + \big(m_2(0)-2m_1^2\big)\,\expo^{-\frac{1}{3}t}$. More generally, we proceed by induction to show that $m_k(t)$ converges exponentially, more precisely $m_k(t)$ is of the form:
  \begin{equation}
    \label{eq:induction_hyp}
    m_k(t) = m_k^*+\mathcal{O}(\expo^{-\frac{k-1}{k+1}t})
  \end{equation}
  with $m_k^*$ the limit value of $m_k(t)$. We first re-write the evolution equation of $m_k(t)$:
  \begin{equation}
    \label{eq:m_k_bis}
    m_k'(t) = -\frac{k-1}{k+1} m_k(t) + P_{k-1}(t)
  \end{equation}
  with $P_{k-1}(t)=\frac{1}{k+1} ∑_{j=1}^{k-1} C_k^jm_j(t) m_{k-j}(t)$. By induction, $P_{k-1}(t)$ has to converge in time. Using variation of constant in \eqref{eq:m_k_bis} gives:
  \begin{equation}
    \label{eq:var_constant}
    m_k(t) = m_k(0) \expo^{-\frac{k-1}{k+1}t} + \expo^{-\frac{k-1}{k+1}t} ∫_{s=0}^t \expo^{\frac{k-1}{k+1}s} P_{k-1}(s)\,\dd s,
  \end{equation}
  which leads to \eqref{eq:induction_hyp}. 
\qed

From the proposition, we observe that the second moment $m_2(t)$ converges exponentially toward the constant $2\,m_1^2$. This behavior could be expected as the equilibrium of the dynamics \eqref{PDE} is given by:
\begin{equation}
  \label{equili}
  q_\infty(x):=\frac{1}{m_1}\expo^{-\frac{x}{m_1}}\mathbbm{1}_{[0,\infty)}(x)
\end{equation}
for which the second moment is equal $2\,m_1^2$.

\begin{remark}
Moment calculations can be useful in the study of classical spatially homogeneous Boltzmann equation, and we refer the readers to \cite{alonso_new_2013} for more information on this regard.
\end{remark}

\subsection{Pairwise distribution}
\label{subsec:2.3}
\begin{figure}[ht]
  \centering
  \includegraphics[width=.8\textwidth]{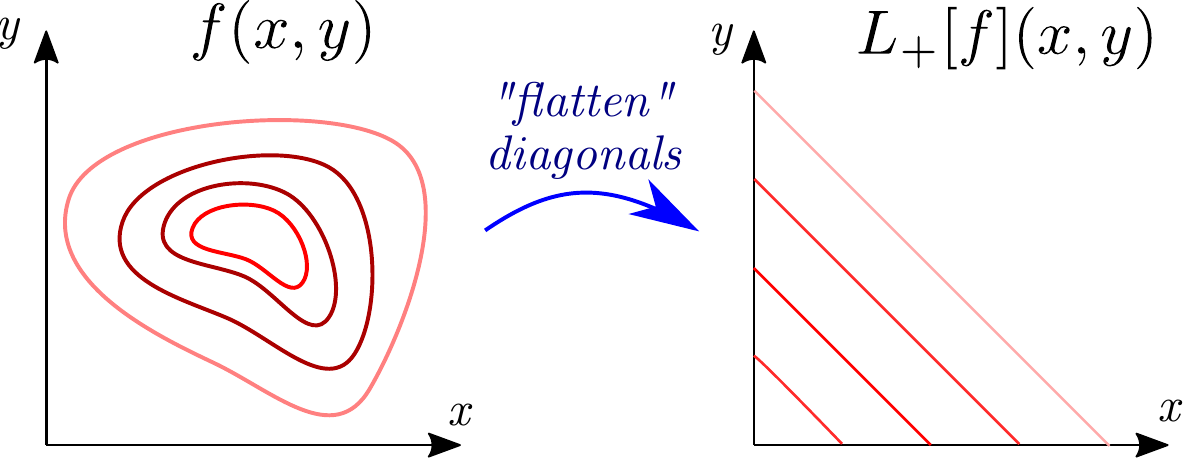}
  \caption{The operator $L_+$ \eqref{eq:L} flattens the distribution $f(x,y)$ over the diagonal lines $x+y=Constant$.}
  \label{fig:flatten}
\end{figure}

Before studying the evolution of the entropy of the solution $q(t,x)$, we make a detour with another formulation of the reshuffling model using a two-particles distribution. Indeed, the jump process $\overline{X}(t)$ \eqref{uniform_limit} is a ``truncated version'' of the following dynamics:
\begin{equation}
  \label{uniform_limit_pair}
  (\overline{X},\overline{Y}) \leadsto \Big(U(\overline{X}\!+\!\overline{Y}) \,,\,(1-U)(\overline{X}\!+\!\overline{Y})\Big)
\end{equation}
where $U\sim \mathrm{Uniform}([0,1])$. Introducing a test function $φ(x,y)$, this dynamics lead to:
\begin{equation}
  \dd 𝔼[φ(\overline{X},\overline{Y})] = 𝔼[φ\Big(U(\overline{X}\!+\!\overline{Y}) \,,\,(1-U)(\overline{X}\!+\!\overline{Y})\Big)-φ(\overline{X},\overline{Y})]\,\dd t.
\end{equation}
We now translate this evolution equation into a PDE.
\begin{proposition}
  Let $f(t,x,y)$ the density distribution of the process $(\overline{X}(t),\overline{Y}(t))$. It satisfies (weakly) the linear evolution equation:
  \begin{equation}
    \label{eq:evo_f}
    ∂_t f = L_+[f] - f
  \end{equation}
  with
  \begin{equation}
    \label{eq:L}
    L_+[f](x,y) = \frac{1}{x+y}∫_{z= 0}^{x+y} f(z,x+y-z) \,\dd z.
  \end{equation}
\end{proposition}
\Proof
The evolution equation \eqref{uniform_limit_pair} gives:
\begin{eqnarray}
  \frac{\dd}{\dd t} ∫_{x,y\geq 0}\!\! f(t,x,y)φ(x,y)\,\dd x\dd y &=&  ∫_{u=0}^1∫_{x,y\geq 0}\!\!\! f(t,x,y)φ\Big(u(x\!+\!y),(1\!-\!u)(x\!+\!y)\Big)\,\dd x\dd y \dd u \nonumber \\
  &&\quad - ∫_{x,y\geq 0} f(t,x,y)φ(x,y)\,\dd x\dd y
\end{eqnarray}
To identify the operator associated with the equation, let's rewrite the ``gain term'' (dropping the dependency in time for simplicity) using two changes of variables:
\begin{eqnarray*}
  && ∫_{u=0}^1∫_{x,y\geq 0}\!\!\! f(x,y)φ\Big(u(x\!+\!y),(1\!-\!u)(x\!+\!y)\Big)\,\dd x\dd y \dd u \hspace{2cm}\\
&& \hspace{2cm} = ∫_{u=0}^1∫_{m\geq 0}∫_{z=0}^m\!\! f(z,m-z)φ\Big(um,(1\!-\!u)m)\Big)\,\dd z\dd m \dd u \\
&& \hspace{2cm} = ∫_{x',y'\geq 0}∫_{z=0}^{x'+y'}\!\! f(z,x'+y'-z)φ(x',y')\,\frac{1}{x'+y'}\dd z\dd x' \dd y'
\end{eqnarray*}
with $(x'=um,y'=(1-u)m)$ leading to $\dd x'\dd y'=m\dd u\dd m$.
\qed
\begin{remark}
  Notice that the operator $L$ \eqref{eq:L} ``flattens'' the distribution $f$ over the diagonals $x+y=\text{Constant}$ and thus minimizing its entropy over each diagonal (see Figure~\ref{fig:flatten}). In particular, the equilibrium for the dynamics are the distribution of the form: $f_*(x,y)= ϕ(x+y)$.
\end{remark}

The linear operator $L_+$ \eqref{eq:L} is linked to the non-linear operator $Q_+$ \eqref{eq:Q_plus}. Indeed, assuming $\overline{X}$ and $\overline{Y}$ are independent, i.e. $f(x,y)=q(x)q(y)$, integrating $L_+[f]$ over the 'extra' variable $y$ gives:
\begin{eqnarray*}
  ∫_{y\geq 0} L_+[f](x,y)\,\dd y &=& ∫_{y\geq 0} \frac{1}{x+y}∫_{z= 0}^{x+y} q(z)q(x+y-z) \,\dd z\,\dd y \\
  &=& ∫_{m=x}^{+∞} \frac{1}{m}∫_{z= 0}^{m} q(z)q(m-z) \,\dd z\,\dd y \quad = Q_+[q](x).
\end{eqnarray*}

\begin{figure}[ht]
  \centering
  \includegraphics[width=.8\textwidth]{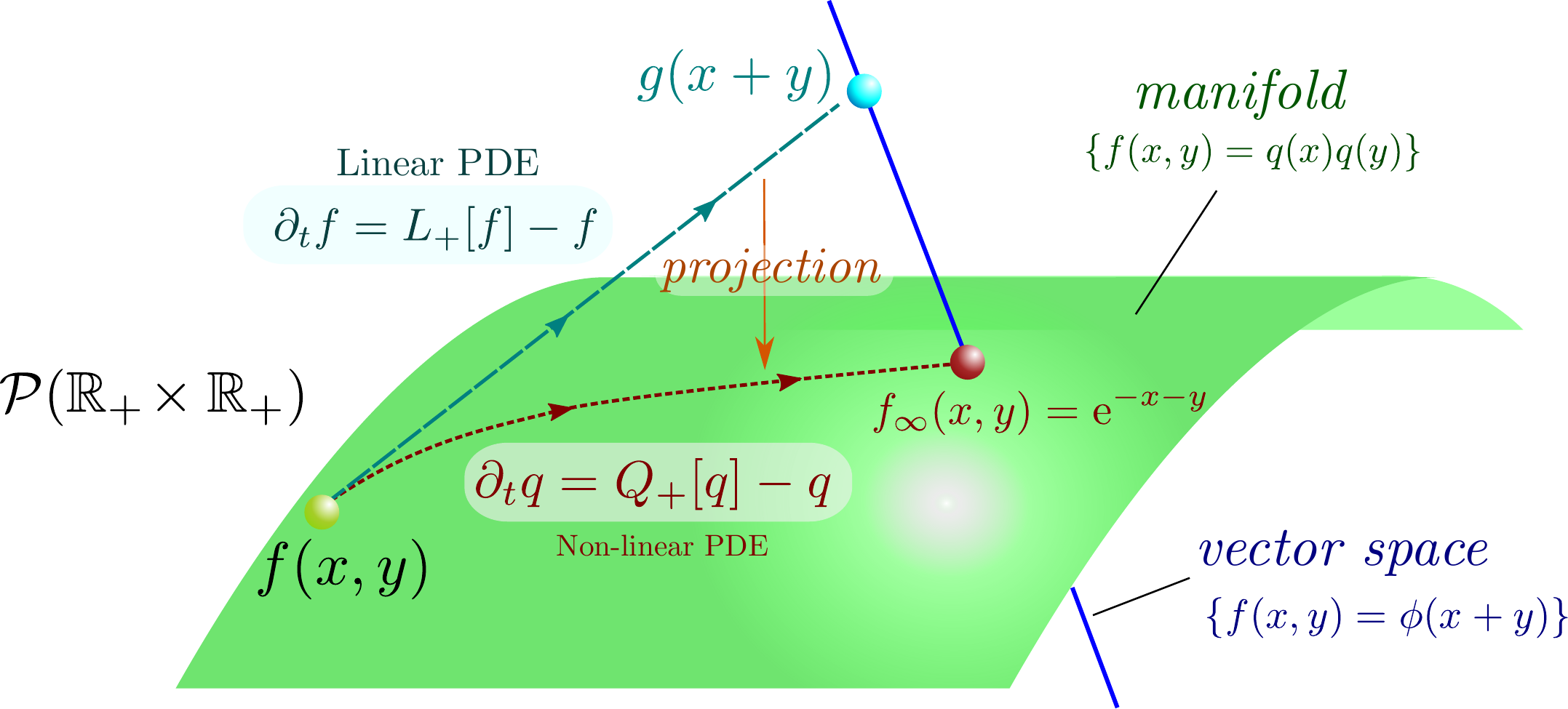}
  \caption{Schematic representation of the evolution of $f(t,x,y)$ and $q(t,x)$. If $f$ belongs to the manifold of independent functions, i.e. $f(t,x,y)=q(t,x)q(t,y)$, then the evolution of its marginal $q$ satisfies {\it locally} the non-linear equation \eqref{PDE}. Notice that the manifold of independent function is not invariant by the flow of the linear PDE. Notice that we have assumed $m_1 = 1$ so that $f_\infty(x,y) := q_\infty(x)q_\infty(y) = \expo^{-x-y}$. Also, the definition of $g$ appears in \eqref{eq:g}.}
  \label{fig:sketch}
\end{figure}

\subsection{Micro-reversibility} 
\label{subsec:2.4}
The evolution equation for $f$ \eqref{eq:evo_f} corresponds to a collisional operator with the kernel:
\begin{equation}
  \label{eq:K}
  K\big(x,y;x',y'\big) = \frac{1}{x+y}δ_{x+y}(x'+y')
\end{equation}
where $δ$ denote the Dirac distribution.
Indeed, writing ${\bf z}=(x,y)$, the equation \eqref{eq:evo_f} could be written:
\begin{equation}
  ∂_t f ({\bf z},t) = ∫_{\tilde{\bf z}\geq 0} K\big(\tilde{\bf z};{\bf z}\big)f(\tilde{\bf z},t)\,\dd \tilde{\bf z} - ∫_{{\bf z}'\geq 0} K\big({\bf z};{\bf z}'\big)f({\bf z},t)\,\dd {\bf z}'
\end{equation}
where ${\bf z}'=(x',y')$ denotes the post-collisional position and $\tilde{\bf z}=(\tilde{x},\tilde{y})$ the pre-collision position.

\begin{remark}
  A more rigorous way to define the kernel $K$ is through a weak formulation using a test function $φ(x,y)$:
\begin{equation}
  "∫_{x',y'\geq 0} K\big(x,y;x',y'\big)φ(x',y')\,\dd x' \dd y'" = \frac{1}{x+y} ∫_{z=0}^{x+y} φ(z,x+y-z) \,\dd z.
\end{equation}
\end{remark}

The collisional kernel $K$ satisfies a micro-reversibility condition, namely:
\begin{equation}
  \label{eq:micro_reversibility}
  K\big({\bf z};{\bf z}'\big)=K\big({\bf z}';{\bf z}\big) \quad \text{for any } {\bf z}\text{ and }{\bf z}'∈ℝ_+×ℝ_+.
\end{equation}

\begin{figure}[ht]
  \centering
  \includegraphics[width=.35\textwidth]{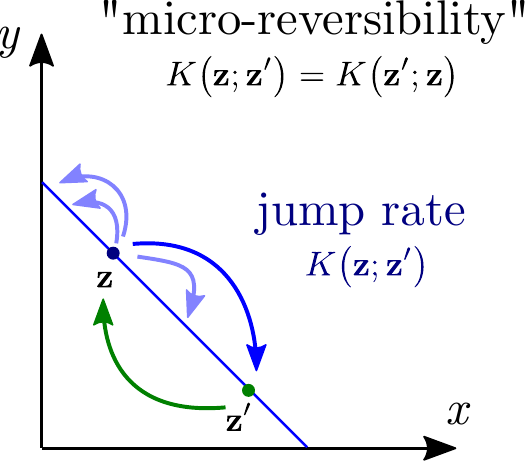}
  \caption{The collisional kernel $K$ \eqref{eq:K} satisfies a micro-reversibility condition.}
  \label{fig:micro_reversibility}
\end{figure}

One has to integrate against a test function $φ$ to make this statement rigorous. As a consequence, we deduce the lemma.
\begin{lemma}
  Let $φ(x,y)$ a (smooth) test function and $f(t,x,y)$ the solution of \eqref{eq:micro_reversibility}. Then:
  \begin{equation}
    \label{eq:f_decay}
        \frac{\dd }{\dd t} ∫_{{\bf z}} f({\bf z},t)φ({\bf z}) \dd {\bf z} = -\frac{1}{2}∫_{{\bf z},{\bf z}'} K({\bf z};{\bf z}') \Big(f({\bf z}',t)-f({\bf z},t)\Big) \Big(φ({\bf z}')-φ({\bf z})\Big) \dd {\bf z} \dd {\bf z}'.
  \end{equation}
\end{lemma}
\Proof
We drop the dependency in time to ease the reading:
\begin{eqnarray*}
  \frac{\dd }{\dd t} ∫_{{\bf z}} f({\bf z})φ({\bf z}) \dd {\bf z} &=& ∫_{\tilde{\bf z},{\bf z}} K(\tilde{\bf z};{\bf z}) f(\tilde{\bf z})φ({\bf z}) \dd \tilde{\bf z} \dd {\bf z} - ∫_{{\bf z},{\bf z}'} K({\bf z};{\bf z}') f({\bf z})φ({\bf z}) \dd \tilde{\bf z} \dd {\bf z} \\
                                                                  &=& ∫_{{\bf z},{\bf z}'} K({\bf z};{\bf z}') f({\bf z})\big(φ({\bf z}')- φ({\bf z})\big) \dd {\bf z} \dd {\bf z}' \\
                                                                  &=& ∫_{{\bf z},{\bf z}'} K({\bf z}';{\bf z}) f({\bf z}')\big(φ({\bf z})- φ({\bf z}')\big) \dd {\bf z} \dd {\bf z}' \\
    &=& ∫_{{\bf z},{\bf z}'} K({\bf z};{\bf z}') f({\bf z}')\big(φ({\bf z})- φ({\bf z}')\big) \dd {\bf z} \dd {\bf z}' \\
  &=& \frac{1}{2}∫_{{\bf z},{\bf z}'} K({\bf z};{\bf z}') \Big(f({\bf z})-f({\bf z}')\Big) \Big(φ({\bf z}')-φ({\bf z})\Big) \dd {\bf z} \dd {\bf z}'.
\end{eqnarray*}
We deduce that both the $L^2$ norm and the entropy of $f(t,x,y)$ decay in time. \qed

\section{Convergence to equilibrium: Wasserstein and linearization}
\label{sec:wasserstein_linearization}
\setcounter{equation}{0}

We carry out an linearization analysis around the exponential equilibrium distribution of the solution of \eqref{PDE} and demonstrate an explicit  rate of convergence under the linearized (weighted $L^2$) setting in section \ref{subsec:3.1}. These arguments are reinforced in section \ref{subsec:3.2} into a local convergence result for the full nonlinear equation. A coupling approach is encapsulated in section \ref{subsec:3.3} in order to show that solution $q(t,x)$ of \eqref{PDE} relaxes to its equilibrium $q_\infty$ exponentially fast in the Wasserstein distance.

\subsection{Linearization around equilibrium}
\label{subsec:3.1}
Now we perform a linearization analysis near the global exponential equilibrium $q_\infty$, in a fashion that is similar to \cite{baranger_explicit_2005}. For this purpose, we define the linear operator $\mathcal{L}$ to be \[\mathcal{L}[h](x):= \int_0^\infty\int_0^\infty \frac{\mathbbm{1}_{[0,k+\ell]}(x)}{k+\ell}q_\infty(k+\ell-x)\Big(h(k)+h(\ell)-h(x)-h(k+\ell-x)\Big)\,\dd k\,\dd \ell.\]
Setting $q = q_\infty(1+\varepsilon h)$ for $0<\varepsilon <<1$, we deduce from \eqref{rewrite} that
\begin{equation}\label{linearized}
\partial_t h(x) = \mathcal{L}[h](x),
\end{equation}
where $h \in L^2(q_\infty)$ is orthogonal to $\mathcal N(\mathcal{L}):=\mathrm{Span}\{1,x\}$ in $L^2(q_\infty)$. For the linearized equation \eqref{linearized}, the natural entropy is the $L^2(q_\infty)$ norm of $h$, and the entropy dissipation is given by
\begin{align*}
  \frac{\dd}{\dd t}\EE :&= \frac{\dd}{\dd t} \frac{1}{2}\|h\|^2_{L^2(q_\infty)} \\
                        &= \int_{\mathbb{R}^3_+ } \frac{\mathbbm{1}_{[0,k+\ell]}(x)}{k+\ell}q_\infty(k)q_\infty(\ell)\big(h(k)+h(\ell)-h(k+\ell-x)-h(x)\big)h(x)\,\dd k \,\dd \ell \, \dd x \\
                        &= -\frac 14\int_{\mathbb{R}^3_+ } \frac{\mathbbm{1}_{[0,k+\ell]}(x)}{k+\ell}q_\infty(k)q_\infty(\ell)\big(h(k+\ell-x)+h(x)-h(k)-h(\ell)\big)^2 \,\dd k \,\dd \ell \,\dd x.
\end{align*}
In particular, it implies that the spectrum of $\mathcal{L}$ in $L^2(q_\infty)$ is non-positive.

\begin{remark}
  It is not hard to show that the linear operator $-\mathcal{L}$ enjoys a self-adjoint property on the space $L^2(q_\infty)$. Thus the existence of a spectral gap $\eta$ is equivalent to \[\forall h \perp \mathcal N(\mathcal{L}),\quad -\langle \mathcal{L}[h],~h\rangle_{L^2(q_\infty)}:= -\int_0^\infty \mathcal{L}[h](x)\, h(x)\, q_\infty(x) \, \dd x \geq \eta \, \|h\|^2_{L^2(q_\infty)}.\]
\end{remark}

\begin{remark}
  Following \cite{grunbaum_linearization_1972}, we give some comments on the space $L^2(q_\infty)$. If $q$ is the unique solution of \eqref{PDE} and we set $q = q_\infty(1+\varepsilon h)$ as before for $h \perp \mathcal N(\mathcal{L})$, then
  \begin{align*}
    \int_0^\infty q\log q \,\dd x &= \int_0^\infty q_\infty(1+\varepsilon h)\, \log(q_\infty(1+\varepsilon h))\, \dd x \\
                                &=\int_0^\infty q_\infty\, \log q_\infty \, \dd x + \varepsilon\int_0^\infty \left(\log\frac{1}{M}-\frac{1}{M}x\right)h \, q_\infty \, \dd x \\
                                &\quad + \int_0^\infty q_\infty(1+\varepsilon h)\left(\varepsilon h-\frac{(\varepsilon h)^2}{2}\pm \cdots \right)\, \dd x\\
                                &=\int_0^\infty q_\infty\, \log q_\infty \, \dd x + \frac{\varepsilon^2}{2}\int_0^\infty h^2\, q_\infty \, \dd x + O(\varepsilon^3),
  \end{align*}
  where we used the fact that $h\perp \mathcal N(\mathcal{L})$. Therefore, we can see that $\|h\|^2_{L^2(q_\infty)} = \int_0^\infty h^2 \, q_\infty \,\dd x$ gives the first-order correction to the expansion of the entropy of $q$ around $q_\infty$.
\end{remark}

We will prove that the linearized entropy $\EE = \frac{1}{2}\|h\|^2_{L^2(q_\infty)}$ decays exponentially fast in time with an explicit sharp decay rate, the essence of which lies in the following lemma.

\begin{lemma}
  Let $m_1=1$ and $\mathcal A:=\{h \in L^2(q_\infty) \mid h\in \mathcal N(\mathcal{L})\}$. Then
  \begin{equation}\label{Key linear}
    \inf\limits_{h\in \mathcal A} \frac{\int_0^\infty h^2(x)\, q_\infty(x)\, \dd x}{\int_0^\infty \frac{\expo^{-z}}{z}\,(\int_0^z h(x)\,\dd x)^2\, \dd z} = 3,
  \end{equation}
  and the infimum in \eqref{Key linear} is attained (up to a non-zero multiplication constant) at $h(x) = \frac{1}{2}\,(x^2-4x+2)$.
\end{lemma}

\Proof
The key ingredient in the proof is the fact that the so-called Laguerre polynomials, defined by \[L_n(x) = \frac{\expo^x}{n!}\frac{\dd^n}{\dd x^n}(\expo^{-x}x^n) = \sum_{k=0}^n \binom{n}{k}\frac{(-1)^k}{k!}x^k,\quad n\geq 0,\] form an orthonormal basis for the weighted $L^2$ space $L^2(q_\infty)$ \cite{abramowitz_handbook_1965}. Thus, for any $h \in L^2(q_\infty)$ which is not identically zero, we can write $h = \sum_{n=0}^\infty \alpha_nL_n$, in which $\alpha_n \in \mathbb R$ for all $n$. Next, notice that the condition $h \in \mathcal A$ implies that $\alpha_0 = \alpha_1 = 0$. Moreover, we have $\int_0^\infty h^2(x)\, q_\infty(x)\, \dd x = \sum_{n=2}^\infty \alpha^2_n$ thanks to the orthonormality of the Laguerre polynomials $\{L_n\}_{n\geq 0}$. To proceed further, we recall that \cite{poularikas_handbook_2018} $L_n(z)-L_{n+1}(z) = \int_0^z L_n(x)\dd x$ and $zL'_n(z) = nL_n(z)-nL_{n-1}(z)$ for all $n\geq 1$, whence
\begin{align*}
  \int_0^\infty \frac{\expo^{-z}}{z}\left(\int_0^z h(x)\dd x\right)^2\dd z &= \int_0^\infty \frac{\expo^{-z}}{z}\left(\sum_{n=2}^\infty \alpha_n(L_n(z)-L_{n+1}(z))\right)^2\dd z \\
                                                                           &= \int_0^\infty \expo^{-z}\left(\sum_{n,m=2}^\infty \alpha_n\alpha_m\left(\frac{L_n(z)-L_{n+1}(z)}{z}\right)(L_m(z)-L_{m+1}(z))\right)^2\dd z \\
                                                                           &= -\sum_{n,m=2}^\infty \frac{\alpha_n\alpha_m}{n+1}\int_0^\infty \expo^{-z}(L_m(z)-L_{m+1}(z))\,\dd L_{n+1}(z) \\
                                                                           &= \sum_{n,m=2}^\infty \frac{\alpha_n\alpha_m}{n+1}\int_0^\infty L_{n+1}(z)\,\dd \big(\expo^{-z}(L_m(z)-L_{m+1}(z))\big) \\
                                                                           &= \sum_{n,m=2}^\infty \frac{\alpha_n\alpha_m}{n+1}\int_0^\infty L_{n+1}(z)\, L_{m+1}(z)\, \expo^{-z}\, \dd z \\
                                                                           &= \sum_{n=2}^\infty \frac{\alpha^2_n}{n+1} \leq \frac{1}{3}\sum_{n=2}^\infty \alpha^2_n.
\end{align*}
Finally, notice that the inequality above will become an equality if and only if $\alpha_n = 0 $ for all $n\geq 3$, or in other words, if and only if $h(x) = L_2(x) = \frac{1}{2}\,(x^2-4x+2)$ up to a non-zero multiplication constant. \qed

We are now in a position to prove the following result.
\begin{theorem}
Assume that $h \in L^2(q_\infty)$ solves the linearized equation \eqref{linearized}, then we have
\begin{equation}\label{eq:expo_decay_for_linearized_eqn}
\|h(t)\|_{L^2(q_\infty)} \leq \|h(0)\|_{L^2(q_\infty)}\expo^{-\frac{1}{3}t}.
\end{equation}
\end{theorem}


\Proof We will only prove the result for $m_1=1$, and the general case follows readily from a change of variable argument. From the discussion above, we already have that
\begin{equation}\label{start}
  \begin{aligned}
    -\frac{\dd}{\dd t} \frac{1}{2}\|h\|^2_{L^2(q_\infty)} &= \int_{\mathbb{R}^3_+ } \frac{\mathbbm{1}_{[0,k+\ell]}(x)}{k+\ell}q_\infty(k)q_\infty(\ell)\cdot \\
    &~~~\left(h(k+\ell-x)+h(x)-h(k)-h(\ell)\right)h(x)\,\dd k \, \dd \ell \, \dd x.
  \end{aligned}
\end{equation}
Thanks to $h\in \mathcal A$, it is not hard to see through a change of variable that \[\int_{\mathbb{R}^3_+ } \frac{\mathbbm{1}_{[0,k+\ell]}(x)}{k+\ell}q_\infty(k)q_\infty(\ell)h(k+\ell-x)h(x)\,\dd k \,\dd \ell \,\dd x = 0.\] Also, a simple calculation yields that
\[\int_{\mathbb{R}^3_+ } \frac{\mathbbm{1}_{[0,k+\ell]}(x)}{k+\ell}q_\infty(k)q_\infty(\ell)h^2(x)\,\dd k \,\dd \ell \,\dd x = \int_0^\infty h^2(x)\, \expo^{-x}\,\dd x \] and \[\int_{\mathbb{R}^3_+ } \frac{\mathbbm{1}_{[0,k+\ell]}(x)}{k+\ell}q_\infty(k)q_\infty(\ell)h(k)h(x)\,\dd k \,\dd \ell \,\dd x = \int_0^\infty \frac{\expo^{-z}}{z}\left(\int_0^z h(x)\,\dd x\right)^2\dd z.\]
Consequently, \eqref{start} reads
\begin{align*}
  -\frac{\dd}{\dd t} \frac{1}{2}\|h\|^2_{L^2(q_\infty)} &= \int_0^\infty h^2(x)\,\expo^{-x}\,
                                                          \dd x - 2\int_0^\infty \frac{\expo^{-z}}{z}\left(\int_0^z h(x)\,\dd x\right)^2\dd z \\
                                                        &\geq \frac{1}{3}\int_0^\infty h^2(x)\,\expo^{-x}\, \dd x = \frac{1}{3}\|h\|^2_{L^2(q_\infty)},
\end{align*}
in which the inequality follows directly from the previous lemma. Thus we can conclude by Gronwall's inequality. \qed
\subsection{Local convergence in $L^2$} 
\label{subsec:3.2}

We now extend the linearization argument from the previous subsection into a local convergence result for the full non-linear equation.
\begin{theorem}
There exists some $\varepsilon>0$ such that if at some time $t\geq 0$,
\[
\int \frac{|q(t,x)-q_\infty(x)|^2}{q_\infty(x)}\,\dd x\leq \varepsilon,
\]
then $q$ converges to $q_\infty$ and for any $\lambda < \frac 13$, there exists $C_\lambda$ such that
\[
\int \frac{|q(t,x)-q_\infty(x)|^2}{q_\infty(x)}\,\dd x\leq C_\lambda\,\expo^{-\lambda\,t}.
\]
\label{localconvergence}
  \end{theorem}
\Proof For a solution $q$, we denote $h(t,x)=(q-q_\infty)/q_\infty$ and calculate
  \[\begin{split}
  &-\frac{\dd }{\dd t} \frac{1}{2}\,\|h\|^2_{L^2(q_\infty)}=-\int h\,\partial_t q =-\int h\,(Q_+[q]-q)\\
  &\quad=-\int h \,q_\infty \,\mathcal{L}[h] -\int h(x)\,q_\infty(x)\,\frac{\mathbbm{1}_{x\leq k+\ell}}{k+\ell}\,q_\infty(k+\ell-x)\,h(k)\,h(\ell)\,\dd x\,\dd k\,\dd \ell.
\end{split}
  \]
Denote \[R(x)=\int \frac{\mathbbm{1}_{x\leq k+\ell}}{k+\ell}\,q_\infty(k+\ell-x)\,h(k)\,h(\ell)\,\dd k\,\dd \ell,\]
and calculate
\begin{align*}
\left|\int h(x)\,q_\infty(x)\,R(x)\,\dd x\right| &\leq \left(\int q_\infty(x)\,\frac{\mathbbm{1}_{x\leq k+\ell}}{k+\ell}\,q_\infty(k+\ell-x)\,h^2(k)\,h^2(\ell)\,\dd x\,\dd k\,\dd \ell\right)^{1/2}\\
&\quad\cdot\left(\int h^2(x)\,q_\infty(x)\,\frac{\mathbbm{1}_{x\leq k+\ell}}{k+\ell}\,q_\infty(k+\ell-x)\,\dd x\,\dd k\,\dd \ell\right)^{1/2}.
\end{align*}
So first of all,
\[
\begin{split}
&\int q_\infty(x)\,\frac{\mathbbm{1}_{x\leq k+\ell}}{k+\ell}\,q_\infty(k+\ell-x)\,h^2(k)\,h^2(\ell)\,\dd x\,\dd k\,\dd \ell\\
&\quad=\int \frac{\mathbbm{1}_{x\leq k+\ell}}{k+\ell}\,q_\infty(k)\,q_\infty(\ell)\,h^2(k)\,h^2(\ell)\,\dd x\,\dd k\,\dd \ell =\|h\|_{L^2(q_\infty)}^4.
\end{split}
\]
On the other hand,
\[\int h^2(x)\,q_\infty(x)\,\frac{\mathbbm{1}_{x\leq k+\ell}}{k+\ell}\,q_\infty(k+\ell-x)\,\dd x\,\dd k\,\dd \ell=\int h^2(x)\,q_\infty(x)\,\dd x = \|h\|_{L^2(q_\infty)}^2.\]
Hence,
\[
\left|\int h(x)\,q_\infty(x)\,R(x)\,\dd x\right|\leq \|h\|_{L^2(q_\infty)}^3.
\]
Coming back to the equation, we have that
\[\begin{split}
  &-\frac{\dd}{\dd t} \frac 12\,\|h\|^2_{L^2(q_\infty)}\geq -\int h(x)\,q_\infty(x)\,\mathcal{L}[h]\,\dd x-\|h\|_{L^2(q_\infty)}^3.
\end{split}
\]
Using the previous calculations on the spectral gap of $\mathcal{L}$, we can conclude that
\[
-\frac{d}{dt} \frac{1}{2}\,\|h\|^2_{L^2(q_\infty)}\geq \frac{1}{3}\,\|h\|_{L^2(q_\infty)}^2-\|h\|_{L^2(q_\infty)}^3,
\]
which finishes the proof with a Gronwall bound.\qed

We can couple this with an interpolation argument to modify the smallness assumption in weighted $L^2$ by using the relative entropy, which leads us to Theorem \ref{thm3} below, whose proof will be deferred to the appendix (as the proof of Theorem \ref{thm3} relies on several a priori estimates established in section \ref{sec:entropy_dissipation}).

\begin{theorem}\label{thm3}
Assume that for some $\lambda_0>\frac 12$, $\sup\limits_x \expo^{\lambda_0 \, x}\,q(0,x)<\infty$. Then there exists some $\delta>0$ such that if at some time $t\geq 0$,
 \[
\int q(t,x)\,\log \frac{q(t,x)}{q_\infty(x)}\,\dd x\leq \delta,
\]
we have that $q$ converges to $q_\infty$ and for any $\lambda < \frac 13$, there exists $C_\lambda$ such that
\[
\int \frac{|q(t,x)-q_\infty(x)|^2}{q_\infty(x)}\,\dd x\leq C_\lambda\,\expo^{-\lambda\,t}.
\]
\label{localconvergenceentropy}
\end{theorem}

\subsection{Coupling and convergence in Wasserstein distance}
\label{subsec:3.3}
In this section we shall employ a coupling argument to demonstrate the convergence of the solution of \eqref{PDE} to the exponential probability density function given by \eqref{equili}. Before we state the main result of this section, we first collect several relevant definitions.

\begin{definition}
For random variables $X$ and $Y$ taking values in $\mathbb R_+$, we write $X \perp Y$ to mean that $X$ and $Y$ are mutually independent. Also, the Wasserstein distance with exponent 2 between two probability density functions (say $f$ and $g$) is defined by
\[W_2(f,g) = \inf\left\{\sqrt{\mathbb E[|X-Y|^2]};~\mathrm{Law}(X)=f,~\mathrm{Law}(Y)=g\right\},\] where the infinimum is taken over all pairs of random variables defined on some probability space $(\Omega,\mathbb P)$ and distributed according to $f$ and $g$, respectively.
\end{definition}

Next, we present a stochastic representation of the evolution equation \eqref{PDE}, which is interesting in its own right.

\begin{proposition}\label{stochastic_representation}
  Assume that $q_t(x) := q(t,x)$ is a solution of \eqref{PDE} with initial condition $q_0(x)$ being a probability density function supported on $\mathbb{R}_+$ with mean $m_1$. Defining $(X_t)_{t\geq 0}$ to be a $\mathbb R_+$-valued continuous-time pure jump process with jumps of the form
  \begin{equation}\label{coupling_nonlinear}
    \begin{array}{ccc}
      X_t & \begin{tikzpicture} \draw [->,decorate,decoration={snake,amplitude=.4mm,segment length=2mm,post length=1mm}]
        (0,0) -- (1,0); \node[above,red] at (0.5,0) {\tiny{rate 1}};\end{tikzpicture} & U(X_t+Y_t),
    \end{array}
  \end{equation}
  where $Y_t$ is a i.i.d. copy of $X_t$, $U \sim \mathrm{Uniform}[0,1]$ is independent of $(X_t)$ and $(Y_t)$, and the jump occurs according to a Poisson clock running at the unit rate. If $\mathrm{Law}(X_0) = q_0$, then $\mathrm{Law}(X_t) = q_t$ for all $t\geq 0$.
\end{proposition}

\Proof
Taking $\varphi$ to be an arbitrary but fixed test function, we have
\begin{equation}\label{testfunc}
  \frac{\dd }{\dd t} \mathbb E[\varphi(X_t)] = \mathbb E[\varphi(U(X_t+Y_t))] - \mathbb E[\varphi(X_t)].
\end{equation}
Denoting $q(t,x)$ as the probability density function of $X_t$, \eqref{testfunc} can be rewritten as
\[\frac{\dd }{\dd t} \int_{\mathbb{R}_+} q(t,x)\varphi(x) \dd x = \int_{\mathbb{R}^2_+}\int_0^1 \varphi(u(k+\ell))q(k,t)q(\ell,t) \dd u \dd k \dd \ell - \int_{\mathbb{R}_+} q(t,x)\varphi(x) \dd x.\] After a simple change of variables, one arrives at
\begin{equation}\label{done}
  \frac{\dd }{\dd t} \int_{\mathbb{R}_+} q(t,x)\varphi(x) \dd x = \int_{\mathbb{R}_+} \left(G[q](x,t) - q(t,x)\right)\varphi(x) \dd x.
\end{equation}
Thus, $q$ must satisfy $\partial_t q = G[q]$ and the proof is completed. \qed

\begin{remark}\label{rem}
  Using a similar reasoning, we can show that if $(\overbar{X}_t)_{t\geq 0}$ is a $\mathbb R_+$-valued continuous-time pure jump process with jumps of the form
  \begin{equation}\label{coupling_limit}
    \begin{array}{ccc}
      \overbar{X}_t & \begin{tikzpicture} \draw [->,decorate,decoration={snake,amplitude=.4mm,segment length=2mm,post length=1mm}]
        (0,0) -- (1,0); \node[above,red] at (0.5,0) {\tiny{rate 1}};\end{tikzpicture} & U(\overbar{X}_t+\overbar{Y}_t),
    \end{array}
  \end{equation}
  where $\overbar{Y}_t$ is a i.i.d. copy of $\overbar{X}_t$, $U \sim \mathrm{Uniform}[0,1]$ is independent of $(\overbar{X}_t)$ and $(\overbar{Y}_t)$, and the jump occurs according to a Poisson clock running at the unit rate. Then $\mathrm{Law}(\overbar{X}_0) = q_\infty$ implies $\mathrm{Law}(\overbar{X}_t) = q_\infty$ for all $t\geq 0$.
\end{remark}

The main result of this section is recorded in the following theorem:

\begin{theorem}\label{Coupling}
  Under the setting of Proposition \ref{stochastic_representation}, we have
  \begin{equation}\label{Wasserstein_nonlinear}
    W_2(q_t,q_\infty) \leq \expo^{-\frac{1}{6}t}W_2(q_0,q_\infty),\quad \forall t\geq 0.
  \end{equation}
\end{theorem}

\Proof
Fixing $t \in \mathbb{R}_+$, we need to couple the two densities $q_t$ and $q_\infty$. Suppose that $(X_t)_{t\geq 0}$ and $(\overbar{X}_t)_{t\geq 0}$ are $\mathbb R_+$-valued continuous-time pure jump processes with jumps of the form \eqref{coupling_nonlinear} and \eqref{coupling_limit}, respectively. We can take $(X_t,Y_t)$ and $(\overbar{X}_t,\overbar{Y}_t)$ as in the statement of Proposition \ref{stochastic_representation} and Remark \ref{rem}, respectively. Meanwhile, we require that $X_t \perp \overbar{Y}_t$, $\overbar{X}_t \perp Y_t$ and $(X_t,\overbar{X}_t) \perp (Y_t,\overbar{Y}_t)$, i.e., several independence assumptions can be imposed along the way when we introduce the coupling. We insist that the same uniform random variable $U$ is used in both \eqref{coupling_nonlinear} and \eqref{coupling_limit}. Moreover, we impose that $\mathrm{Law}(X_0)=q_0$ and $\mathrm{Law}(\overbar{X}_0)=q_\infty$. As a consequence of the previous proposition and remark, $q_t = \mathrm{Law}(X_t)$ and $\mathrm{Law}(\overbar{X}_t)=q_\infty$ for all $t\geq 0$, whence $\mathbb E[\overbar{X}_t]=\mathbb E[\overbar{Y}_t]=m_1$ and $\mathbb E(\overbar{X}^2_t)=\mathbb E(\overbar{Y}^2_t)= 2\,m^2_1$, $\forall t\geq 0$. Also, we have that $\mathbb E[X_t]=\mathbb E[Y_t]=m_1$ for all $t\geq 0$. Thanks to the aforementioned coupling, we then have
\begin{align*}
  \frac{\dd}{\dd t}\mathbb E[(X_t-\overbar{X}_t)^2] &= \mathbb E[\big(U(X_t+Y_t-\overbar{X}_t-\overbar{Y}_t)\big)^2 - (X_t-\overbar{X}_t)^2] \\
                                                    &= \frac 13\Big(\mathbb E[(X_t-\overbar{X}_t)^2]+\mathbb E[(Y_t-\overbar{Y}_t)^2]+2\mathbb E[(X_t-\overbar{X}_t)(Y_t-\overbar{Y}_t)]\Big) \\
                                                    &\qquad \qquad - \mathbb E[(X_t-\overbar{X}_t)^2] \\
                                                    &= \frac 23\mathbb E[(X_t-\overbar{X}_t)^2] + \frac 23\mathbb E[X_t-\overbar{X}_t]\cdot\mathbb E[Y_t-\overbar{Y}_t] - \mathbb E[(X_t-\overbar{X}_t)^2]\\
                                                    &= -\frac 13\mathbb E[(X_t-\overbar{X}_t)^2].
\end{align*}
Now we pick $\overbar{X_0}$ with law $q_\infty$ so that $W^2_2(q,q_\infty) = \mathbb E[(X_0-\overbar{X}_0)^2]$, and a routine application of Gronwall's inequality yields \eqref{Wasserstein_nonlinear}. \qed

\section{Entropy dissipation}
\label{sec:entropy_dissipation}
\setcounter{equation}{0}

We state our main result, Theorem \ref{EEP}, in section \ref{subsec:4.1} so that readers know exactly what is at stake. We will present various expressions of the entropy and entropy dissipation associated to the solution $q(t,x)$ of \eqref{PDE}, along with a discussion of the strategy of the proof of Theorem \ref{EEP} in section \ref{subsec:4.2}. A sequence of auxiliary lemmas and corollaries are recorded in section \ref{subsec:4.3} and \ref{subsec:4.4}. Finally, a full proof of Theorem \ref{EEP}, built upon all of the preparatory work from \ref{subsec:4.1} to \ref{subsec:4.4}, is shown in \ref{subsec:4.5}.

\subsection{Main result}
\label{subsec:4.1}
For the integro-differential equation \eqref{PDE}, a common strategy  \cite{bassetti_explicit_2010,during_kinetic_2008,matthes_steady_2008} is to use the Laplace transform or Fourier transform of \eqref{PDE} to prove the exponential decay of solution of \eqref{PDE} to $q_\infty(x)$ in some Fourier metric. However, little analysis of \eqref{PDE} has been carried out without resorting to Laplace or Fourier transform. In particular, we would like to show the dissipation of relative entropy, i.e., $\mathrm{D_{KL}}(q(\cdot,t)~||~ q_\infty)$, along solution trajectories:
\begin{equation}\label{entropy_dis}
  \frac{\dd}{\dd t} \int_0^\infty q\,\log \frac{q}{q_\infty}\,\dd x = \frac{\dd}{\dd t} \int_0^\infty q\,\log q \,\dd x \leq 0.
\end{equation}
It is reasonable to expect the validity of \eqref{entropy_dis} as the exponential probability density $q_\infty$ maximizes the negative entropy $-\int_0^\infty p\,\log p \,\dd x$ among all continuous probability density functions supported on $[0,\infty)$ with prescribed mean.

The following proposition together with its proof should be a reminiscent of the calculations carried out for a standard Boltzmann equation arising from the kinetic theory of (dilute) gases \cite{villani_review_2002}.

\begin{proposition}\label{weak_formulation}
  Let $\varphi(x)$ be a (continuous) test function on $\mathbb{R}_+$ and assume that $q$ is a smooth solution of \eqref{PDE}, then we have
  \begin{align*}
    \frac{\dd}{\dd t} \int_0^\infty q(t,x)\varphi(x) \dd x = -\frac 14\int_{\mathbb{R}^3_+}& \frac{\mathbbm{1}_{[0,k+\ell]}(x)}{k+\ell}\big(q(k+\ell-x)q(x)-q(k)q(\ell)\big)\cdot\\ &\big(\varphi(k+\ell-x)+\varphi(x)-\varphi(k)-\varphi(\ell) \big)\dd k \dd \ell \dd x.
  \end{align*}
  Moreover, inserting $\varphi = \log q$ and employing the fact that total mass is conserved (i.e., $m'_0(t) = 0$ for all $t\geq 0$), we obtain the dissipation of relative entropy:
  \begin{equation*}
    \frac{\dd}{\dd t} \int_0^\infty q(t,x)\log q(t,x)\dd x = -\frac 14 D[q],
  \end{equation*}
   where \begin{equation}\label{D}
   D[q] := \int_{\mathbb{R}^3_+} \frac{\mathbbm{1}_{[0,k+\ell]}(x)}{k+\ell}\big(q(k+\ell-x)q(x)-q(k)q(\ell)\big)\,
   \log \frac{q(k+\ell-x)q(x)}{q(k)q(\ell)}\,\dd k \,\dd \ell \,\dd x \geq 0.
   \end{equation}
\end{proposition}

\Proof
We notice that the PDE \eqref{PDE} can be rewritten as
\begin{equation}\label{rewrite}
  \partial_t q(x) = \int_0^\infty\int_0^\infty \frac{\mathbbm{1}_{[0,k+\ell]}(x)}{k+\ell}\Big(q(k)q(\ell)-q(x)q(k+\ell-x)\Big)\dd k \dd \ell
\end{equation}
(thanks to Proposition \ref{conservation}). Omitting the time variable for simplicity, we deduce that
\begin{align*}
  \frac{\dd}{\dd t} \int_0^\infty q(x)\varphi(x)\dd x &= \int_{\mathbb{R}^3_+ } \frac{\mathbbm{1}_{[0,k+\ell]}(x)}{k+\ell}\Big(q(k)q(\ell)-q(x)q(k+\ell-x)\Big)\varphi(x)\dd k \dd \ell \dd x \\
                                                      &= \int_{\mathbb{R}^3_+ } \frac{\mathbbm{1}_{[0,k+\ell]}(x)}{k+\ell}q(k)q(\ell)\Big(\varphi(x)-\varphi(\ell)\Big) \dd k \dd \ell \dd x\\
                                                      &= \int_{\mathbb{R}^3_+ } \frac{\mathbbm{1}_{[0,k+\ell]}(x)}{k+\ell}q(k)q(\ell)\Big(\varphi(k+\ell-x)-\varphi(k)\Big) \dd k \dd \ell \dd x\\
                                                      &= \frac 12\int_{\mathbb{R}^3_+ } \frac{\mathbbm{1}_{[0,k+\ell]}(x)}{k+\ell}q(k)q(\ell) \\
                                                      &\qquad\cdot \Big(\varphi(k+\ell-x)+\varphi(x)-\varphi(k)-\varphi(\ell)\Big) \dd k \dd \ell \dd x\\
                                                      &= -\frac 14\int_{\mathbb{R}^3_+ } \frac{\mathbbm{1}_{[0,k+\ell]}(x)}{k+\ell}\Big(q(k+\ell-x)q(x)-q(k)q(\ell)\Big)\\
                                                      &\qquad\cdot \Big(\varphi(k+\ell-x)+\varphi(x)-\varphi(k)-\varphi(\ell)\Big) \dd k \dd \ell \dd x.
\end{align*}
\qed

\begin{remark}
  The dissipation of the relative entropy can also be seen via an alternative perspective. Indeed, we fix $t \geq 0$ and assume that $X_1(t)$ and $X_2(t)$ are i.i.d $\mathbb{R}_+$-valued random variable with its probability density function given by $q(t,x)$, and we define $(Z_1,Z_2) = (U(X_1+X_2),(1-U)(X_1+X_2))$ with $U \sim \mathrm{Uniform}[0,1]$  being independent of $X_1$ and $X_2$. Then we deduce from the PDE \eqref{PDE} and Lemma \ref{lem1} that
  \begin{equation}\label{entro_decay}
    \begin{aligned}
      2\frac{\dd}{\dd t}\mathrm{D_{KL}}(q ~||~ q_\infty) &= \HH((Z_1,Z_2),(X_1,X_2)) - \HH((X_1,X_2)) \\
      &\leq \HH((Z_1,Z_2)) - \HH((X_1,X_2)),
    \end{aligned}
  \end{equation}
  where $\HH(X,Y) := \int_{\mathbb{R}} \rho_X(x)\log \rho_Y(x) \dd x$ represents the cross entropy from $Y$ to $X$, if the laws of $X$ and $Y$ are given by $\rho_X$ and $\rho_Y$. It can be shown \cite{apenko_monotonic_2013} that the joint entropy of $(Z_1,Z_2)$ is always no more than the joint entropy of $(X_1,X_2)$, whence the rightmost side of \eqref{entro_decay} is non-positive.
\end{remark}


\begin{corollary}
The  exponential distribution $q_\infty$ defined in \eqref{equili} is the only (smooth) equilibrium solution of the PDE \eqref{PDE}.
\end{corollary}

\Proof
By Proposition \ref{weak_formulation}, we see that \[q_\infty(x)q_\infty(k+\ell-x) = q_\infty(k)q_\infty(\ell) \quad \text{for all $k,\ell,x\geq 0$ such that $k+\ell \geq x$}.\] Since $\int_0^\infty q_\infty(x) \dd x= 1$ and $\int_0^\infty xq_\infty(x) \dd x = m_1$, $q_\infty$ must be the exponential probability density provided by \eqref{equili}. \qed

We will prove that $q \xrightarrow[]{t \to \infty} q_\infty$ polynomial in time. Without of loss generality, throughout the argument to be presented below we will set $m_1 = 1$, i.e., $q_\infty(x) = \expo^{-x}$ for $x \geq 0$. Our main result is stated as follows:

\begin{theorem}\label{EEP}
  Under the assumptions of Lemma \ref{Laplace} below, we have for some constant $C,\,\theta>0$ and for any $t\geq C\, \log (1\slash D)$ that
  \begin{equation}
    \label{eq:main_result}
    \int_{x=0}^{+∞} q(x,t)\,\log\frac{q(x,t)}{\expo^{-x}}\,\dd x\leq C\,D^\theta.
  \end{equation}
\end{theorem}
To our best knowledge, Theorem \ref{EEP} is the first entropy-entropy dissipation inequality established for the uniform reshuffling dynamics.

\subsection{Basic expressions of the entropy-entropy dissipation}
\label{subsec:4.2}
Let us start by looking at the strong convergence of the pairwise distribution, which is essentially trivial. Indeed, we recall the linear PDE \eqref{eq:evo_f}, which reads
\begin{equation*}
\partial_t f=L_+[f]-f,
\end{equation*}
where
\[L_+[f](x,y)=\frac{1}{x+y}\,\int_{z=0}^{x+y} f(z,x+y-z)\,\dd z.\]
Then denoting
\begin{equation}\label{eq:g}
g(t,\lambda)=\frac{1}{\lambda}\,\int_0^\lambda f(t,z,\lambda-z)\dd z,
\end{equation}
we can rewrite \eqref{eq:evo_f} as $\partial_t f(t,x,y)=g(t,x+y)-f(t,x,y)$, whence
\begin{align*}
\partial_t g(t,\lambda) &=\frac{1}{\lambda}\,\int_{z=0}^{\lambda} \partial_t f(t,z,\lambda-z)\,\dd z\\
&=\frac{1}{\lambda}\,\int_0^\lambda \left(g(t,\lambda)-f(t,z,\lambda-z)\right)\,\dd z = 0.
\end{align*}
Hence $g(t,\lambda)=g(0,\lambda)$ and trivially
\begin{equation}\label{convergencef}
|f(t,x,y)-g(0,x+y)|\leq \expo^{-t}.
\end{equation}
Unfortunately this cannot be used to show the convergence on the actual equation for $q(t,x)$ because the two models are not equivalent: If $q(t,x)$ solves \eqref{PDE}, which is nonlinear, then in general $f(t,x,y)=q(t,x)\,q(t,y)$ does not solve \eqref{eq:evo_f}. The one exception is when $q(t,x)$ is some exponential.

This can also be seen from the fact that in the argument above $f$ does not necessarily converge to an exponential but to whatever $g(t=0)$ was. The rate of convergence is also too fast as the second moment of $q$ converges much slower for example.

We will still find some of the structure above in the entropy dissipation for $q$ but that is one reason why the entropy dissipation is not easy to handle. In particular, the entropy dissipation will vanish whenever $f(x,y)=g(x+y)$ which seems to create some degeneracy.

Next, we can rewrite the dissipation term in a manner that will make the connection with the exponential more apparent. We define for simplicity $f(x,y)=q(x)\,q(y)$, and as before
\[g(\lambda)=\frac{1}{\lambda}\,\int_0^\lambda f(z,\lambda-z)\,\dd z=\frac{1}{\lambda}\,\int_0^\lambda q(z)\,q(\lambda-z)\,\dd z.\]
Finally, we also define
\[h(x)=\int_{\R_+} g(x+y)\,\dd y.\]
We remark here that $h$ coincides with the collision gain operator $Q_+[q]$ defined via \eqref{eq:Q_plus}. With these definitions, we have
\begin{lemma}
One has that
\[D = 2\,\int_{\R_+^2} q(x)\,q(y)\,\log \frac{q(x)\,q(y)}{g(x+y)}\,\dd x\,\dd y+  2\,\int_{\R_+^2} g(x+y)\,\log \frac{g(x+y)}{q(x)\,q(y)}\,\dd x\,\dd y,\]
or as well that
\[
\begin{split}
D &=2\,\int_{\R_+^2} q(x)\,q(y)\,\log \frac{q(x)\,q(y)}{g(x+y)}\,\dd x\,\dd y+  2\,\int_{\R_+^2} g(x+y)\,\log \frac{g(x+y)}{h(x)\,h(y)}\,\dd x\,\dd y\\
  &~~+4\,\int_{\R_+} h(x)\,\log \frac{h(x)}{q(x)}\,\dd x.
\end{split}
\]
\end{lemma}
Formally this forces $g(x+y)$ to be close to $f(x,y)$ so this is a very similar term to the one that we had found when looking at equation \eqref{eq:evo_f}. It is some sort of degeneracy because it does not directly force $f$ to be close to $\expo^{-x-y}$ so we will have to resolve it. Of course since $f(x,y)=q(x)\,q(y)$, $f(x,y)=g(x+y)$ forces $q$ to be some exponential and therefore this should be possible.

\Proof We can first simply rewrite
\[
D=\int_{\R_+^3} \frac{\mathbbm{1}_{y+z\geq x}}{y+z}\,(f(y+z-x,x)-f(y,z))\,\log \frac{f(y+z-x,x)}{f(y,z)}\,\dd x\,\dd y\,\dd z.
\]
Observe that by swapping $x$ and $z$
\begin{align*}
&\int_{\R_+^3} \frac{\mathbbm{1}_{y+z\geq x}}{y+z}\,(f(y+z-x,x)-f(y,z))\,\log f(y+z-x,x)\\
&=\int_{\R_+^3} \frac{\mathbbm{1}_{y+x\geq z}}{y+x}\,(f(y+x-z,z)-f(y,x))\,\log f(y+x-z,z).
\end{align*}
Changing variable $y\to y'=y+x-z$, we get that
\begin{align*}
&\int_{\R_+^3} \frac{\mathbbm{1}_{y+z\geq x}}{y+z}\,(f(y+z-x,x)-f(y,z))\,\log f(y+z-x,x)\\
&=\int_{\R_+^3} \frac{\mathbbm{1}_{y'+z\geq x}}{y'+z}\,(f(y',z)-f(y'+z-x,x))\,\log f(y',z).
\end{align*}
Hence \[D=2\,\int_{\R_+^3} \frac{\mathbbm{1}_{y+z\geq x}}{y+z}\,(f(y,z)-f(y+z-x,x))\,\log f(y,z).\]
In other words,
\[D=2\int_{\R_+^2} f(y,z)\,\log f(y,z)\,dy\,dz-2\,\int_{\R_+^2} g(y+z)\,\log f(y,z)\,\dd y\,\dd z.\]
Now, we observe that \[\int_{\R_+^2} f(y,z)\,\log g(y+z)\, \dd y \, \dd z = \int_{\R_+^2} g(y+z)\log g(y+z)\, \dd y \, \dd z.\] Indeed, a change of variable $y = x - w$ and $z = w$ yields \[\int_{\R_+^2} g(y+z)\log g(y+z)\, \dd y \, \dd z = \int_{\R_+} x\, g(x) \, \log g(x) \, \dd x.\] By the same change of variables, we also have \[\int_{\R_+^2} f(y,z)\,\log g(y+z)\, \dd y \, \dd z = \int_{\R_+} \log g(x) \, \int_0^x f(x-w,w) \, \dd w \dd x = \int_{\R_+} x\, g(x) \, \log g(x) \, \dd x.\]
Hence
\[\frac{D}{2}=\int_{\R_+^2} f(y,z)\,\log \frac{f(y,z)}{g(y+z)}\,\dd y\,\dd z+\int_{\R_+^2} g(y+z)\,\log \frac{g(y+z)}{f(y,z)}\,\dd y\,\dd z.\]
Finally as $f(y,z)=q(y)\,q(z)$, we may also notice that
\begin{align*}
\int_{\R_+^2} g(y+z)\,\log \frac{g(y+z)}{f(y,z)}\,\dd y\,\dd z &=\int_{\R_+^2} g(y+z)\,\log g(y+z)\,\dd y\,\dd z\\
&~~~-2\,\int_{\R_+^2} g(y+z)\,\log q(y)\,\dd y\,\dd z\\
&=\int_{\R_+^2} g(y+z)\,\log g(y+z)\,\dd y\,\dd z-2\,\int_{\R_+} h(y)\,\log q(y)\,\dd y.
\end{align*}
So we also have that
\begin{align*}
\int_{\R_+^2} g(y+z)\,\log \frac{g(y+z)}{f(y,z)}\,\dd y\,\dd z &=\int_{\R_+^2} g(y+z)\,\log \frac{g(y+z)}{h(y)\,h(z)}\,\dd y\,\dd z\\
&~~~+2\,\int_{\R_+} h(y)\,\log \frac{h(y)}{q(y)}\,\dd y,
\end{align*}
concluding the estimate. \qed

Next, we intend to collect here some various bounds stemming from the dissipation term, the essence of those bounds lies in the following lemma.

\begin{lemma}\label{Miscellaneous_estimates}
We have that \[\int q(x)\,\log\frac{q(x)}{H(x)}\,dx\leq \int \varphi(y)\,q(x)\,q(y)\,\log \frac{q(x)\,q(y)}{g(x+y)}\,\dd x\,\dd y,\]
in which \[H(x)=\int g(x+y)\,\varphi(y)\,\dd y, \] for any $\varphi \geq 0$ such that $\int \varphi\,q\, \dd x=1$.
\end{lemma}

\Proof  Indeed, as $\log$ is concave,
\begin{align*}
\int q(x)\,\varphi(y)\,q(y)\,\log \frac{g(x+y)}{q(x)\,q(y)}\,\dd x\,\dd y &\leq \int q(x) \log\left(\int \frac{g(x+y)}{q(x)} \,\varphi(y)\,dy\right) \,\dd x\\
&=\int q(x)\,\log \frac{H(x)}{q(x)}\,\dd x,
\end{align*}
and the proof is completed. \qed

As a consequence of this lemma, inserting $\phi(x)=1$ and then $\phi(x)=x$, we then deduce that
\begin{align*}
&\int q(x)\,\log\frac{q(x)}{h(x)}\,\dd x\leq \int q(x)\,q(y)\,\log \frac{q(x)\,q(y)}{g(x+y)}\,\dd x\,\dd y,\\
&\int q(x)\,\log\frac{q(x)}{m(x)}\,\dd x\leq \int x\,q(x)\,q(y)\,\log \frac{q(x)\,q(y)}{g(x+y)}\,\dd x\,\dd y,
\end{align*}
where \[m(x)=\int g(x+y) \, y \, \dd y = \int_x^\infty g(z) \, (z-x) \, \dd z = \int_x^\infty \,\int_y^\infty g(z) \, \dd z\, \dd y = \int_x^\infty h(y)\,\dd y.\]

\begin{remark}
We also note that $m(0)=1$ (since $\int h\,\dd x=\int q\,\dd x=1$) and so
\[\int h\,\log m\,\dd x=-\int m'\,\log m \, \dd x= -\int h\,\dd x= -\int x\, h(x)\,\dd x = -1,\]
by virtue of the fact that $\int x\, h(x)\,\dd x=\int x\,q(x)\,\dd x=1$.
Thus, \[\int h\,\log \frac{h}{m}\,\dd x=\int h\,\log \frac{h}{\expo^{-x}}\,\dd x.\]

This leads to a possible strategy: Control $\int h\,\log \frac{h}{m}$ in terms of $\int q\,\log \frac{q}{h}$, $\int h\,\log \frac{h}{q}$ and $\int q\,\log \frac{q}{m}$. Then control $\int q\,\log \frac{q}{\expo^{-x}}$ by the previous quantities and $\int h\,\log \frac{h}{m}$. We can then estimate $\int x\,q(x)\,q(y)\,\log \frac{q(x)\,q(y)}{g(x+y)}\,\dd x\,\dd y$ via $\int q(x)\,q(y)\,\log \frac{q(x)\,q(y)}{g(x+y)}\,\dd x\,\dd y$ and some control on the decay of $q$ at infinity. So in the end this would lead to some kind of bounds on $\int q\,\log \frac{q}{\expo^{-x}}$ in terms of the dissipation term. We illustrate the strategy in Figure \ref{fig:sketch2_bis}.

\begin{figure}[ht]
  \centering
  \includegraphics[width=.75\textwidth]{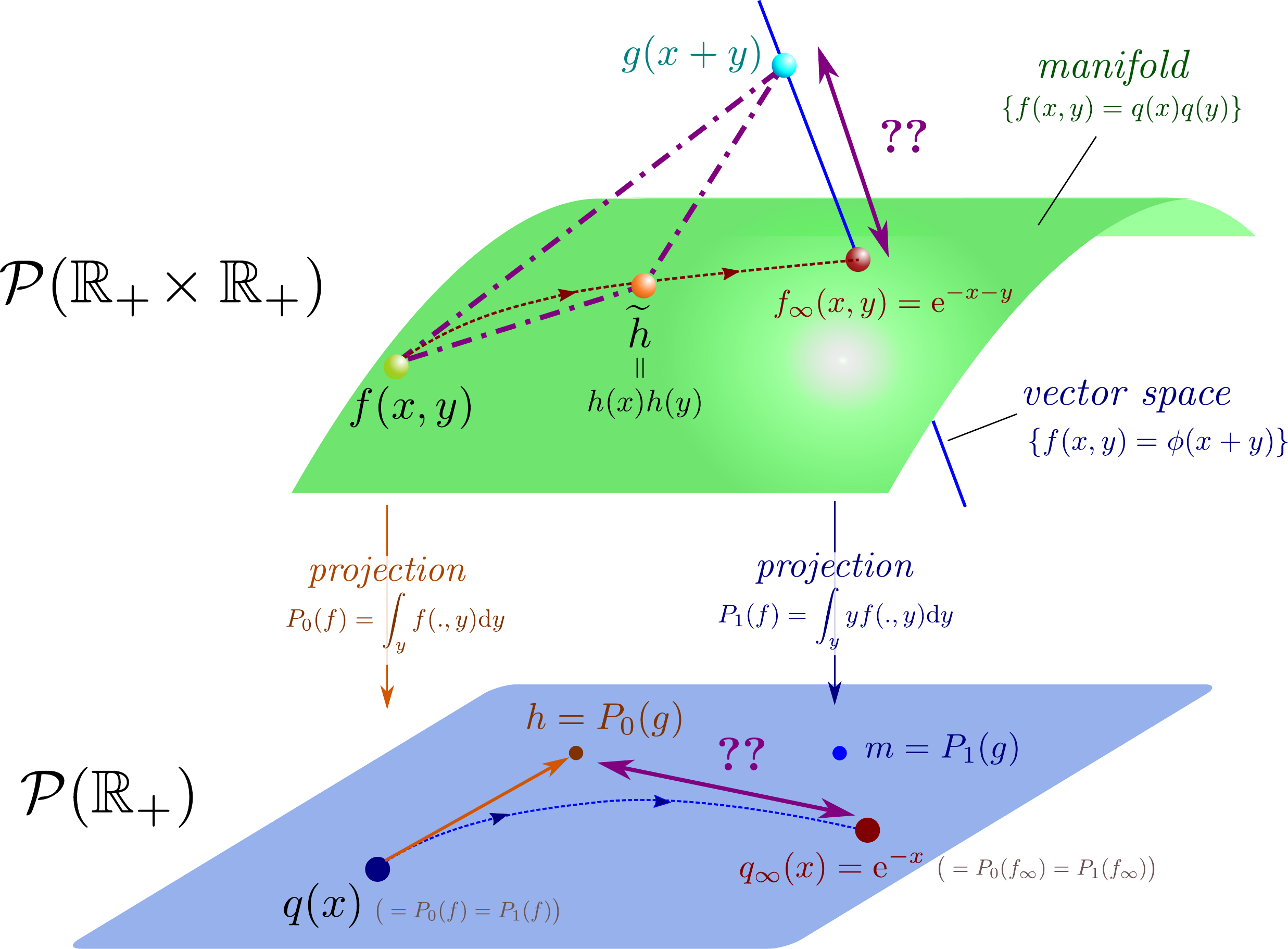}
  \caption{To measure the decay of the relative entropy $\int q\,\log \frac{q}{\expo^{-x}}$, we have to control the term $\int h\,\log \frac{h}{\expo^{-x}}$ or similarly the term $\int g\,\log \frac{g}{\expo^{-x-y}}$ (represented in purple). Indeed, the dissipation term $D$ already provides a control over the 'triangle' of relative entropies $∫ f \log \frac{f}{g}$, $∫ g \log \frac{g}{\widetilde{h}}$ and $∫  \widetilde{h}\log \frac{\widetilde{h}}{g}$ with $\widetilde{h}(x,y)=h(x)h(y)$.}
  \label{fig:sketch2_bis}
\end{figure}

However normally it is not possible to switch relative entropy estimates. Indeed, it is not so hard to find examples of non-negative functions $\varphi,\;\phi,\;\psi$ with total mass $1$ such that
\[\int \varphi\,\log \frac{\varphi}{\psi}=\infty,\]
while \[\int \phi\,\log \frac{\phi}{\psi}+\int \phi\,\log \frac{\phi}{\varphi}+\int \varphi\,\log \frac{\varphi}{\phi}<\infty.\]
Therefore this strategy is not obvious to implement. It should work nicely if we had a control like $\expo^{-x} \slash C\leq q(x)\leq C\,\expo^{-x}$ but the general case is certainly trickier. What saves us is the key observation that here $h$ and $m$ are actually very nice functions in all cases. For example, $m$ and $h$ are monotone decreasing so bounded from above and bounded from below on any finite interval (from the propagation of moments on $q$). This gives us some hope when implementing the aforementioned machinery. We emphasize here that our entropy-entropy dissipation argument draws inspiration from earlier works on Becker-D\"{o}ring equations and coagulation models \cite{jabin_rate_2003,canizo_trend_2017}.
\end{remark}

\subsection{Switching relative entropies}
\label{subsec:4.3}
We note that the relative entropy behaves in the following manner
\begin{lemma}\label{label1}
For any two $\mu,\nu\in \mathcal{P}(\R_+)$ and for any $C\geq 2$, then
\begin{equation}\label{orderofentropy}
\begin{aligned}
&\frac{1}{2\,C}\,\int_{\nu/C\leq \mu\leq C\,\nu} \frac{(\mu-\nu)^2}{\nu}+\frac{1}{8}\,\int_{\mu\leq \nu/C} \nu+\frac{1}{4}\,\int_{\mu\geq C\,\nu} \mu\,\log \frac{\mu}{\nu}\\
&\leq \int \mu\,\log\frac{\mu}{\nu}\\
&\leq \frac{C}{2}\,\int_{\nu/C\leq \mu\leq C\,\nu} \frac{(\mu-\nu)^2}{\nu}+\int_{\mu\leq \nu/C} \nu+\,\int_{\mu\geq C\,\nu} \mu\,\log \frac{\mu}{\nu}.
\end{aligned}
\end{equation}
\end{lemma}

\Proof We observe that
\[\int \mu\,\log\frac{\mu}{\nu}=\int \nu\,\left(\frac{\mu}{\nu}\,\log\frac{\mu}{\nu}+1-\frac{\mu}{\nu}\right).\]
On the other hand, around $1$, the function $\phi(x)=x\,\log x+1-x$ satisfies that $\phi(x)\leq (x-1)^2/2$ for $x\geq 1$ and $\phi(x)\leq \frac{C}{2}\,(x-1)^2$ for $1/C\leq x\leq 1$. On the other hand $\phi(x)\geq (x-1)^2/2C$ for $1\leq x\leq C$ and $\phi(x)\geq (x-1)^2/2$ for $1/C\leq x\leq 1$.
Furthermore $\phi$ lies between  $1/8$ and $1$ when $x\leq 1/2$ and larger than $\frac{x}{4}\,\log x$ for $x\geq 2$. \qed

\begin{remark}
One can also rewrite a little bit the statement of Lemma \ref{label1} so that we do not need to impose that $\mu$ and $\nu$ are probability measures.
\end{remark}

This allows us to ``switch'' relative entropies between two measures that are comparable.
\begin{corollary}\label{changeentropy}
There exists a constant $C>0$ such that if  $\mu_1,\,\mu_2,\,\nu\in \mathcal{P}(\R_+)$ with $\lambda^{-1}\,\mu_1\leq \mu_2\leq \lambda\,\mu_1$ and $\lambda \geq \expo$, then
\[\int \mu_1 \,\log\frac{\mu_1}{\nu}\leq C\,\lambda^3\,\int \mu_2 \,\log\frac{\mu_2}{\nu}+\lambda^3\int \mu_2\,\log\frac{\mu_2}{\mu_1}.\]
\end{corollary}

\Proof Apply Lemma \ref{label1} with $C=2\,\lambda$ first on $\mu_1$ and $\nu$ to find
\[\int \mu_1\,\log\frac{\mu_1}{\nu} \leq \lambda\,\int_{\frac{\nu}{2\lambda}\leq \mu_1\leq 2\lambda\,\nu} \frac{(\mu_1-\nu)^2}{\nu}+\int_{\mu_1\leq \frac{\nu}{2\lambda}} \nu +\,\int_{\mu_1\geq 2\,\lambda\,\nu} \mu_1\,\log \frac{\mu_1}{\nu}.\]
Thanks to Lemma \ref{label1} again, we have \[\int_{\mu_1\leq \frac{\nu}{2\lambda}} \nu \leq 8\,\int \mu_2 \,\log\frac{\mu_2}{\nu}.\]
Now if $\mu_1\leq \frac{\nu}{2\lambda}$ then $\mu_2\leq \frac{\nu}{2}$. Similarly if $\mu_1\geq 2\,\lambda\,\nu$ then $\mu_2\geq 2\,\nu$ and moreover
\[\mu_1\,\log \frac{\mu_1}{\nu}\leq \lambda\,\mu_2\,\log \frac{\lambda\,\mu_2}{\nu}\leq 3\lambda\,\log\lambda\,\mu_2\,\log \frac{\mu_2}{\nu}.\]
Conversely if $\frac{\nu}{2\lambda}\leq \mu_1\leq 2\lambda\,\nu$ then $\frac{\nu}{2\lambda^2}\leq \mu_2\leq 2\lambda^2\,\nu$, and
\[\frac{(\mu_1-\nu)^2}{\nu}\leq 2\,\left(\frac{(\mu_2-\nu)^2}{\nu}+\frac{(\mu_1-\mu_2)^2}{\nu}\right)\leq 2\,\frac{(\mu_2-\nu)^2}{\nu}+4\,\lambda\,\frac{(\mu_1-\mu_2)^2}{\mu_1}.\]
Hence
\begin{align*}
\int_{\frac{\nu}{2\lambda}\leq \mu_1\leq 2\lambda\,\nu} \frac{(\mu_1-\nu)^2}{\nu}&\leq 2\,\int_{\frac{\nu}{2\lambda^2}\leq \mu_2\leq 2\lambda^2\,\nu} \frac{(\mu_2-\nu)^2}{\nu}\\
&~~~+4\lambda\,\int_{\frac{\mu_1}{\lambda}\leq \mu_2\leq \lambda\,\mu_1} \frac{(\mu_1-\mu_2)^2}{\mu_1}.
\end{align*}
Note that by Lemma \ref{label1} applied with $C=\lambda$, we have that
\[
\int_{\frac{\mu_1}{\lambda}\leq \mu_2\leq \lambda\,\mu_1} \frac{(\mu_1-\mu_2)^2}{\mu_1}\leq 2\,\lambda\,\int \mu_2\,\log\frac{\mu_2}{\mu_1}.
\]
Also, Lemma \ref{label1} applied with $C=2\,\lambda^2$ gives rise to
\[\int_{\frac{\nu}{2\lambda^2}\leq \mu_2\leq 2\lambda^2\,\nu} \frac{(\mu_2-\nu)^2}{\nu}\leq 4\,\lambda^2\,\int \mu_2\,\log \frac{\mu_2}{\nu}.\]
Assembling these estimates, the proof is completed. \qed

\subsection{Additional a priori estimates}
\label{subsec:4.4}
This leads us to try to compare $q$ and $h$. We first observe that we can get easy upper bounds.
\begin{lemma}\label{Laplace}
Assume that for some $0<\lambda_0<1$, $\int \expo^{\lambda_0x}\,q(t=0,x)\,\dd x<\infty$. Then we have that
\[\sup_t \int \expo^{\lambda_0\,x}\,q(t,x)\,\dd x<\infty.\]
\end{lemma}

\Proof We use a Laplace transform by defining
\[F(t,\lambda)=\int \expo^{\lambda\,x}\,q(t,x)\,\dd x,\]
and note that
\[
\partial_t F=\int_{\R^2_+} \frac{\expo^{\lambda\,(y+z)}-1}{\lambda\,(y+z)}\,q(y)\,q(z)\,\dd y\,\dd z-F=\frac{1}{\lambda}\, \int_0^\lambda (F(\mu))^2\,\dd \mu-F.
\]
It is useful to remark right away that the stationary solution to this equation satisfies that $F^2=\partial_\lambda (\lambda\,F)$ which has solutions of the form $\frac{1}{1-C\,\lambda}$. Those do blow-up but only for $\lambda$ large enough. As a matter of fact since $\partial_\lambda F|_{\lambda=0}=1$, we can see that we should even have $C=1$.
For this reason, denote now $G=(1-C\,\lambda)\,F$ with some $C < \frac{1}{\lambda}$ such that $G(t=0,\lambda)\leq 1$ on $[0,\ \lambda_0]$. We first show that $\sup\limits_{\lambda \in [0,\lambda_0]} G(t,\lambda) \leq 1$ for all $t\geq 0$. Indeed, let $\lambda(t)$ be such that $\sup\limits_{\lambda \in [0,\lambda_0]} G(t,\lambda) = G(t,\lambda(t))$, then \[\partial_t \sup_{\lambda \in [0,\lambda_0]} G(t,\lambda) \leq \partial_t G(t,\lambda(t)) ,\] this is because $\partial_\lambda G(t,\lambda(t)) = 0$ if $\lambda(t) < \lambda_0$, while if $\lambda(t) = \lambda_0$ then $\partial_\lambda G(t,\lambda(t)) \leq 0$ and $\lambda'(t) \leq 0$, leading to the same inequality. Now since
\begin{equation}\label{eqn:G}
\partial_t G=(\lambda^{-1}-C)\,\int_0^\lambda \frac{(G(\mu))^2}{(1-C\, \mu)^2}\,\dd \mu-G,
\end{equation}
together with $\int_0^\lambda \frac{d\mu}{(1-C\,\mu)^2}=\frac{\lambda}{1-C\,\lambda}$, we deduce that
\[\partial_t \sup_{\lambda \in [0,\lambda_0]} G(t,\lambda) \leq \left(\sup_{\lambda \in [0,\lambda_0]} G(t,\lambda)\right)^2 - \sup_{\lambda \in [0,\lambda_0]} G(t,\lambda),\] which yields via the maximum principle that $\sup\limits_{\lambda \in [0,\lambda_0]} G(t,\lambda) \leq 1$. Now thanks to \eqref{eqn:G} again and the elementary observation that $\partial_t \sup\limits_{\lambda \in [0,\lambda_0]} G(t,\lambda) \leq \sup\limits_{\lambda \in [0,\lambda_0]} \partial_t G(t,\lambda)$, we arrive at
\[\partial_t \sup_{\lambda\in [0,\ \lambda_0]} G(\lambda)\leq 0,\]
which immediately proves the desired upper bound. \qed

\begin{remark}
We believe it is possible to prove the exponential convergence of the Laplace transform $F(t,\lambda)$ to $1\slash (1-\lambda)$ over $\lambda\in [0,\lambda_0)$. However, this is not strictly better though than having the exponential convergence in some weak Wasserstein norm plus the control of the exponential moments that is given above, so we did not try too much in this direction.
\end{remark}

Out of Lemma \ref{Laplace}, we may deduce pointwise bounds on $q$ and $h$, for this purpose, we need the following preparatory result.

\begin{lemma}\label{unifbound:h(t,0)}
We have that \[\sup_{t\geq 0} h(t,0) < \infty,\] i.e., $h(t,0)$ is uniformly bounded in time.
\end{lemma}

\Proof To show $h(t,0)$ is uniformly bounded in time, we write
\begin{align*}
h(t,0) &= \int_{\R^2_+} \frac{q(y)\, q(z)}{y+z} \, \dd y \, \dd z = 2\, \int\int_{y \leq z} \frac{q(y)\, q(z)}{y+z} \, \dd y \, \dd z \\
&\leq 2\, \int\int_{y\leq z} \frac{q(y)\, q(z)}{z} \, \dd y \,\dd z \\
&\leq 2\, \sup_{y\leq r} \, \int_{z \geq r} \frac{r\, q(z)}{z} \, \dd z + 2\, \sup_{y\leq r} \, \int\int_{y\leq z,\, z\leq r} \frac{q(z)}{z} \, \dd y \, \dd z + \frac{2}{r}.
\end{align*}
We know that there exists some $r$ uniformly in time such that \[\int\int_{y\leq z,\, z\leq r} \frac{q(z)}{z} \, \dd y \, \dd z = \int_{z\leq r} q(z) \, \dd z \leq \frac 18.\]  Moreover, for this $r$ we also have $\int_{z\geq r} \frac{r\, q(z)}{z} \, \dd z \leq \frac 18$. Thus,  \[h(t,0) \leq \frac{1}{2} \, \sup_{x\leq r} q(x) + \frac{2}{r}.\] Now we recall the equation for $q$ to find that for any $x \leq r$, \[\partial_ t q(t,x) \leq h(t,0)-q(t,x) \leq \frac 12 \, \sup_{x \leq r} q(t,x) + \frac{2}{r} - q(t,x),\] so if $x_*$ is such that $q(t,x_*) =  \sup\limits_{x \leq r} q(t,x)$, then \[\partial_ t q(t,x_*) ≤ \frac{2}{r} - \frac 12 \, q(t,x_*).\] By Gronwall's inequality, we deduce that $\sup\limits_{x ≤ r} q(t,x) \leq \frac{4}{r}$, which allows us to finish the proof. \qed

\begin{corollary}\label{exponentialpointwise}
Assume that for some $0<\lambda_0<1$, $\int \expo^{\lambda_0 x}\,q(0,x)\,\dd x<\infty$, then we have that
\[
\begin{split}
&\sup_t \int \expo^{\lambda_0 x}\,h(t,x)\,\dd x < \infty, \, \sup_{t,x} \expo^{\lambda_0 x}\,h(t,x)<\infty,\\
& q(t,x)\leq C\,\expo^{-\lambda_0\,x}+q(0,x)\,\expo^{-t}\ \mbox{for some}\ C>0.
\end{split}
\]
\end{corollary}

\Proof  The first bound follows from the definition of $h$. Indeed, as $h=Q_+[q]$, we have
\[\begin{split}
\int \expo^{\lambda_0\, x}\,h(t,x)\,\dd x&=\int \frac{\expo^{\lambda_0\,(y+z)}-1}{\lambda_0\,(y+z)}\,q(y)\,q(z)\,\dd y\,\dd z\\
&\leq \int \expo^{\lambda_0\,(y+z)}\,q(y)\,q(z)\,\dd y\,\dd z < \infty.
\end{split}\]
Next we observe that $h$ is decreasing in $x$, so for any $x \geq 0$
\[\begin{split}
\int_0^\infty \expo^{\lambda_0\,y}\,h(t,y)\,\dd y \geq \int_0^x \expo^{\lambda_0\,y}\,h(t,y)\,\dd y &\geq h(t,x)\,\int_0^x \expo^{\lambda_0\,y}\,\dd y\\
&=h(t,x)\,\frac{\expo^{\lambda_0\,x}-1}{\lambda_0}.
\end{split}\]
Since $h(t,x)\leq h(t,0)$ is uniformly bounded in time, this shows the second point.

\noindent Finally we recall the equation for $q$, which reads $\partial_t q= h-q$, so we may rewrite \eqref{PDE} as
\begin{equation}\label{inteqq}
q(t,x)=q(0,x)\,\expo^{-t}+\int_0^t h(s,x)\,\expo^{-(t-s)}\,\dd s.
\end{equation}
Moreover, notice that
\[\begin{split}
\expo^{\lambda_0\,x}\,\int_0^t h(s,x)\,\expo^{-(t-s)}\,\dd s &\leq \sup_s \left(\expo^{\lambda_0\,x}\,h(s,x)\right)\,\int_0^t \expo^{-(t-s)}\,\dd s\\
&\leq \sup_s \left(\expo^{\lambda_0\,x}\,h(s,x)\right).
\end{split}.\]
Combining these estimates with \eqref{inteqq} ends the proof. \qed

We now turn to lower bounds on $q$ and hence $h$. We start with a lower bound on $q$ in terms of $h$.
\begin{lemma}\label{boundqht-1}
There exists $C$ such that for any $t\geq 1$,
\begin{equation}
q(t,x)\geq \frac{1}{C}\,h(t-1,x).
\end{equation}
\end{lemma}

\Proof We note from the equation \eqref{PDE} that
\[\partial_t h(t,x)=2\,\int_{x}^\infty \frac{1}{\lambda}\,\int_0^\lambda h(t,z)\,q(t,\lambda-z)\,\dd z\,\dd \lambda -2\,h(t,x).\]
Therefore
\[\partial_t h(t,x)\geq -2\,h(t,x), \]
and we have that for any $s\leq t$ that
\[\partial_t q(t,x)\geq \expo^{-(t-s)}\,h(s,x)-q(t,x),\]
leading for example to the claimed result \[q(t,x)\geq \frac{h(t-1,x)}{C}\] with $C = \frac{\expo^2}{\expo-1}$, thereby completing the proof.\qed

\noindent Unfortunately, this is not enough to give us a bound between $q$ and $h$ which would solve everything. Instead, we can first deduce a bound near the origin.
\begin{lemma}\label{boundorigin}
There exists a constant $C$ such that
\begin{equation}
\inf_{t\geq 1}\inf_{x\in [0,\ 2]} h(t,x)\geq \frac{1}{C},\qquad \inf_{t\geq 2}\inf_{x\in [0,\ 2]} q(t,x)\geq \frac{1}{C}.
\end{equation}
\end{lemma}

\Proof For any $x\leq 2$, we have that
\[\begin{split}
h(t,x)&=\int \frac{\mathbbm{1}_{x\leq y+z}}{y+z}\,q(t,y)\,q(t,z)\,\dd y\,\dd z\\
&\geq \int_{y,z\geq 1} \frac{1}{(y+1)\,(z+1)}\,q(y)\,q(z)\,\dd y\,\dd z=\left(\int_1^\infty \frac{q(y)}{1+y}\,\dd y\right)^2.
\end{split}\]
By Cauchy-Schwartz, we have that
\[\begin{split}
\int_1^\infty q(y)\,\dd y&\leq \left(\int_1^\infty \frac{q(y)}{1+y}\,\dd y\right)^{1/2}\,\left(\int_1^\infty (1+y)\,q(y)\,\dd y\right)^{1/2}\\
&\leq \left(\int_1^\infty \frac{q(y)}{1+y}\,\dd y\right)^{1/2}\,\left(\int_0^\infty (1+y)\,q(y)\,\dd y\right)^{1/2}\\
&=\sqrt{2}\,\left(\int_1^\infty \frac{q(y)}{1+y}\,\dd y\right)^{1/2}.
\end{split}\]
On the other hand the convergence of all moments of $q$ shows that there exists $C$ such that for all $t\geq 1$,
\[\int_1^\infty q(y)\,dy\geq \frac{1}{C}.\]
Therefore there exists $C$ such that $h(t,x)\geq \frac{1}{C}$ whenever $x\leq 2$ and $t\geq 1$.
Finally, we deduce the second result from Lemma \ref{boundqht-1}. \qed

We combine this with the following doubling type of argument.
\begin{lemma}\label{lemdoubling0}
There exists a constant $C$ such that for any $x$ and $t \geq 1$, there holds
\[q(t,x)\geq \frac{x}{C}\,\left(\inf_{s\in [t-1,t]} \inf_{y\in [x/2,\;3x/4]} q(s,y)\right)^2.\]
\end{lemma}

\Proof This is a simple consequence of a lower bound on $h$. Indeed, we have
\[\begin{split}
h(t,x)&=\int \frac{\mathbbm{1}_{x\leq y+z}}{y+z}\,q(t,y)\,q(t,z)\,\dd y\,\dd z\\
&\geq \frac{2}{3\,x}\,\int_{y,\,z\in [x/2,\;3x/4]} q(y)\,q(z)\,\dd y\,\dd z.
\end{split}\]
Therefore, \[h(t,x)\geq \frac{x}{24}\,\left(\inf_{y\in [x/2,\;3x/4]} q(t,y)\right)^2.\]
We can again conclude by virtue of Lemma \ref{boundqht-1}. \qed

\begin{lemma}\label{lemdoubling}
There exists a constant $C$ such that for any $t \geq 2$ and $x \geq 2$, we have
\[q(t,x)\geq  \int_{y\geq x} \frac{q(t-1,y)}{C\,y} \, \dd y.\]
\end{lemma}

\Proof This is again a consequence of a lower bound on $h$. Indeed,
\[\begin{split}
h(t,x)&=\int_{\R^2_+} \frac{\mathbbm{1}_{x\leq y+z}}{y+z}\,q(t,y)\,q(t,z)\,\dd y\,\dd z\\
&\geq \int_{y\leq x}\int_{z\geq x} q(y)\,\frac{q(z)}{2\,z}\,\dd y\,\dd z.
\end{split}\]
Thus, by the lower bound on $q$ on $[0,2]$ (thanks to Lemma \ref{boundorigin}), we arrive at
\[h(t,x)\geq \int_{y\geq x} \frac{q(t,y)}{C\,y} \, \dd y.\]
Using Lemma \ref{boundqht-1}, we can again conclude. \qed

\noindent Owing to Lemma \ref{lemdoubling}, we immediately deduce that
\begin{corollary}\label{explowerbound}
There exists some $C>0$ such that for any $x \geq 2$ and any $t\geq \max(C\, x,1)$
\[h(t,x)\geq \frac{\expo^{-C\,x}}{C},\quad q(t,x)\geq \frac{\expo^{-C\,x}}{C}.\]
\end{corollary}

\begin{figure}[ht]
\centering
\includegraphics[scale=0.9]{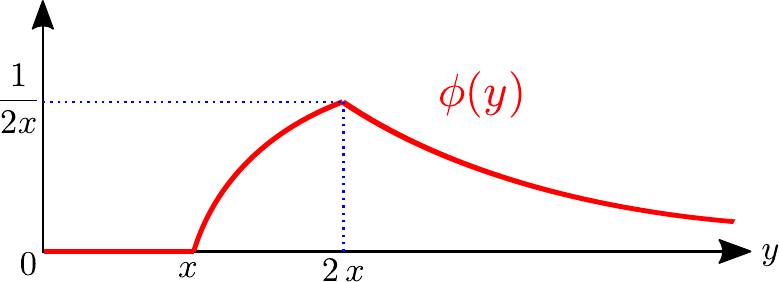}
\caption{The function $\phi$ used in the proof of Corollary \ref{explowerbound}. Notice that $\phi(y) \leq \frac{1}{y}$ for all $y>0$.}
\label{fig:phi}
\end{figure}

\Proof Define $\phi(y)=\frac{(y/x-1)_+}{y}$ for $y\leq 2\,x$ and $\phi=1/y$ if $y\geq 2\,x$ (see Figure \ref{fig:phi}).

\noindent Note that $\phi$ is Lipschitz with
\[\|\nabla \phi\|_{L^\infty}\leq \frac{1}{x^2}.\]
Hence \[x^2\,\int \phi(y)\,q(y)\,\dd y\geq x^2\, \int \phi(y)\,\expo^{-y}\,\dd y-W_1(q,\expo^{-x}),\]
in which $W_1(q,\expo^{-x})$ represents the Wasserstein distance (with exponent 1) between $q$ and $\expo^{-x}$. Thanks to the exponentially fast in time of the convergence $W_1(q,\expo^{-x}) \to 0$, which is a simple consequence of Theorem \ref{Coupling}, we deduce that
\[\int_{y\geq x} \frac{q(y)}{y}\,\dd y\geq \int \phi(y)\,\expo^{-y}\,\dd y-\frac{C}{x^2}\,\expo^{-t/6}.\]
Note that
\[\int \phi(y)\,\expo^{-y}\,\dd y\geq \int_{y\geq 2x} \frac{\expo^{-y}}{y}\,\dd y\geq \frac{\expo^{-3x}}{3x}\,\int_{2x\leq y\leq 3x} \dd y=\frac{\expo^{-3x}}{3}.\]
Therefore from Lemma \ref{lemdoubling}, we can conclude provided that $\frac{C}{x^2}\,\expo^{-t/6}\leq \frac{\expo^{-3x}}{6}$. \qed

\subsection{Proof of the main result}
\label{subsec:4.5}
Armed with all the previous estimates, we can finally present the proof of Theorem \ref{EEP}.

\Proof \!\!\!{\bf of theorem \ref{EEP}}. We note from our earlier estimates that
\[\begin{split}
&\int q\,\log\frac{q}{m}\,\dd x\leq \int x\,q(x)\,q(y)\,\log\frac{q(x)\,q(y)}{g(x+y)}\,\dd x\,\dd y\\
&=\int \left(x\,q(x)\,q(y)\,\log\frac{q(x)\,q(y)}{g(x+y)} +x\,g(x+y)-x\,q(x)\,q(y)\right)\,\dd x\,\dd y.
\end{split}\]
For some $K>0$, we can separate the integral into those $x\leq K$, for which
\[\int_{x\leq K} \left(x\,q(x)\,q(y)\,\log\frac{q(x)\,q(y)}{g(x+y)} +x\,g(x+y)-x\,q(x)\,q(y)\right)\,\dd x\,\dd y\leq K\,D.\]
On the other hand, denoting $\phi(x)=x\,\log x+1-x$, which is a non-negative convex function on $\R_+$ and satisfies $\phi(x) \leq C\, x$ for some constant $C$ if $x$ is bounded, we deduce for any $\lambda \in (0,\lambda_0)$ that
\begin{align*}
&\int_{x\geq K} \left(x\,q(x)\,q(y)\,\log\frac{q(x)\,q(y)}{g(x+y)} +x\,g(x+y)-x\,q(x)\,q(y)\right)\,\dd x\,\dd y\\
&~~~\leq \frac{1}{\lambda}\,\int_{x\geq K} g(x+y)\,\phi\circ\phi\left(\frac{q(x)\,q(y)}{g(x+y)}\right)\,\dd x\,\dd y +\frac{1}{\lambda}\, \int_{x\geq K} \expo^{\lambda \,x}\,g(x+y)\,\dd x\,\dd y,
\end{align*}
where the inequality follows from the Fenchel's inequality $x\, y \leq \phi(x) + \phi^*(y)$, in which $\phi^*$ denotes the Legendre convex conjugate of $\phi$ (and one can check that $\phi^*(y) = \expo^y - 1 \leq \expo^y$ and also $\left(\frac{\phi}{\lambda}\right)^*(y) \leq \frac{\expo^{\lambda \, y}}{\lambda}$).

\noindent We can immediately note that $\phi\circ\phi\leq x\log x$ for large $x$. Thus from Corollary \ref{exponentialpointwise}, we have that
\[\begin{split}
&\int_{x\geq K} \left(x\,q(x)\,q(y)\,\log\frac{q(x)\,q(y)}{g(x+y)} +x\,g(x+y)-x\,q(x)\,q(y)\right)\,\dd x\,\dd y\\
&~~~\leq \frac{D}{\lambda}+\frac{C}{\lambda}\,\expo^{-(\lambda_0-\lambda)\,K}.
\end{split}\]
Combining both estimates gives rise to \[\int q\,\log\frac{q}{m}\leq (K+1)\,\frac{D}{\lambda}+\frac{C}{\lambda}\,\expo^{-(\lambda_0-\lambda)\,K},\]
and optimizing in $K$ leads to
\begin{equation}\label{qmvsD}
\int q\,\log\frac{q}{m}\leq C\,D\,\log\frac{1}{D}.
\end{equation}
The next step is to change this to $\int h\,\log \frac{h}{m}$. We decompose again
\[\int h\,\log \frac{h}{m}=\int_{x\leq K} (h\,\log \frac{h}{m}+m-h)+\int_{x\geq K} (h\,\log \frac{h}{m}+m-h).\]
We note that since $h=-\partial_x m$,
\[\int_{x\geq K} h\,\log m=-\int_{x\geq K} \partial_x m\,\log m=m(K)\,\log m(K)-m(K).\]
Applying Corollary \ref{exponentialpointwise} again, this shows that for some constant $C$, we have that
\begin{equation}\label{changeentropylargex}
\int_{x\geq K} (h\,\log \frac{h}{m}+m-h)\leq C\,\expo^{-K/C}.
\end{equation}
From Corollary \ref{exponentialpointwise} and Corollary \ref{explowerbound}, we note that on $x\leq K$ there holds $\expo^{-C\,K}\leq \frac{q}{h}\leq \expo^{C\,K}$, at least provided that $t\geq C\, x$. As we will see soon, we will choose $K$ logarithmic in $1\slash D$ which gives the assumption appearing in the statement of Theorem \ref{EEP}.

Now in the region $x\leq K$, we can use Lemma \ref{label1} in exactly the same manner as what we did in Corollary \ref{changeentropy}, which yields that
\begin{align*}
\int (h\,\log \frac{h}{m}+m-h) &\leq C\,\expo^{C\,K}\,\int q\,\log \frac{q}{m}\,+C\,\expo^{-K/C} \\
&\leq C\,\expo^{C\,K}\, D\,\log\frac{1}{D} \, + C\,\expo^{-K/C}.
\end{align*}
Optimizing in $K$, we find that for some $\theta>0$ (but $\theta<1$ unfortunately),
\begin{equation}\label{changeqtoh}
\int h\,\log \frac{h}{m}\leq C\,\left(D\,\log\frac{1}{D}\right)^\theta.
\end{equation}
Now we recall that, as a simple consequence of Lemma \ref{Miscellaneous_estimates}, we have
\begin{equation}\label{hmtoexp}
\int h\,\log \frac{h}{m}=\int h\,\log \frac{h}{\expo^{-x}}.
\end{equation}
Therefore we now want to change back from $h$ to $q$. This is the same process and leads to
\begin{equation}\label{changehtoq}
\int q\,\log \frac{q}{\expo^{-x}}\leq C\,\left(\int h\,\log \frac{h}{\expo^{-x}}\right)^{\theta}.
\end{equation}
To finish the proof, we just need to combine \eqref{changehtoq} with \eqref{hmtoexp} and \eqref{changeqtoh}.
\qed

\begin{figure}[!htb]
\centering
\includegraphics[scale=0.8]{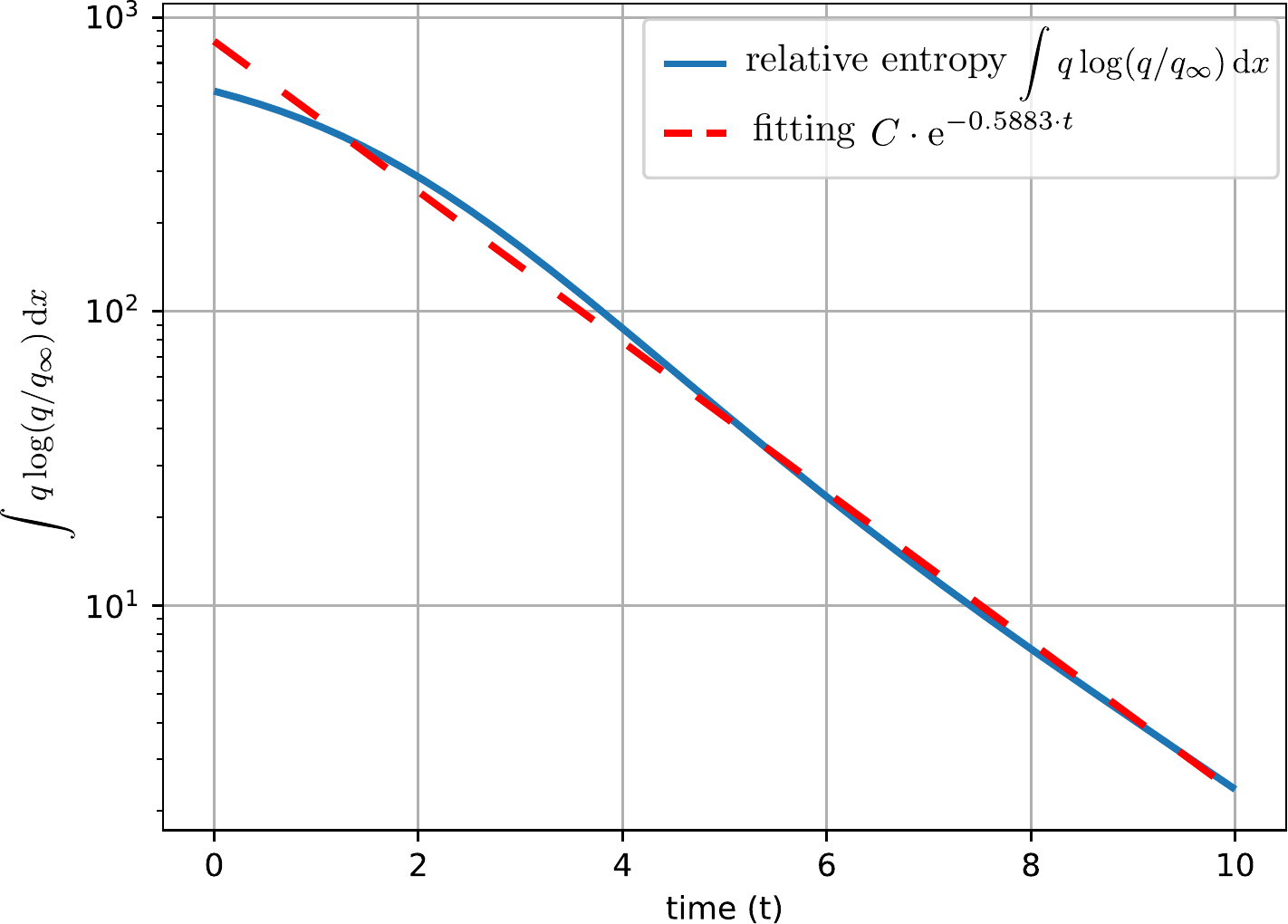}
\caption{Simulation of the relative entropy from $q$ to $q_\infty$ after $t = 10$ in the semilogy scale. We employed forward Euler method with time step-size $\Delta t = 0.05$, space step-size $\Delta x = 0.01$, and a ``random'' initial condition $q(t=0,x)$ having mean value $m_1 = 5$ for the numerical simulation of \eqref{PDE}. This experiment suggests that the relaxation of $\int q\,\log(q\slash q_\infty)\,\dd x$ might be exponentially fast in time, instead of polynomially fast in time as guaranteed by Theorem \ref{EEP}.}
\label{entropy_decay}  
\end{figure}

We end this section with a numerical experiment demonstrating the entropic convergence of $q$ to $q_\infty$, see Figure \ref{entropy_decay}.


\section{Propagation of chaos}
\label{sec:propagation of chaos}
\setcounter{equation}{0}

We give the statement of the propagation of chaos, Theorem \ref{PoC} in \ref{subsec:5.1}. A technical lemma that will be employed in the proof of Theorem \ref{PoC} is displayed in \ref{subsec:5.2}. We reveal the full proof of Theorem \ref{PoC} in \ref{subsec:5.3}.

\subsection{Statement of propagation of chaos}
\label{subsec:5.1}
In this section, we try to adapt the martingale-based techniques developed in \cite{merle_cutoff_2019,hermon_cutoff_2020} to justify the propagation of chaos \cite{sznitman_topics_1991}. For this purpose, we equip the space $\mathcal{P}(\mathbb{R}_+)$ with the Wasserstein distance with exponent $1$, which is defined via \[W_1(\mu,\nu) = \sup\limits_{\|\nabla \varphi\|_\infty \leq 1} \langle \mu-\nu, \varphi\rangle\] for $\mu,\nu \in \mathcal{P}(\mathbb{R}_+)$. We will also need the following version of It\^o's formula.

\begin{lemma}\label{Ito}
Consider an inhomogeneous Poisson process $\NN_t$ with intensity $\lambda(t)$, and a random variable $Y(t)$ left-continuous and adapted to the filtration $\mathcal{F}_t$ generated by $\NN_t$. We define the compound jump process $Z(t)$ and $M(t)$ its associated compensated martingale by:
\begin{equation}\label{Z}
\dd Z(t) = Y(t)\,\dd \NN_t,\quad M(t) = Z(t) - Z(0) - \int_0^t \tilde Y(s)\,\lambda(s)\, \dd s,
\end{equation}
where $\tilde Y$ is any other left-continuous and adapted process.
It\^o's lemma then implies that for any $C^1$ function $\Phi$,
\begin{equation}\label{Itoformula}
\dd \mathbb{E}[\Phi(M(t))] = \mathbb{E}\left[\Phi\left(M(t-)+Y(t)\right)- \Phi(M(t-))\right]\lambda(t)\,\dd t - \mathbb{E}[\nabla\Phi'(M(t))\cdot \tilde Y(t)\,\lambda(t)]\dd t.
\end{equation}
\end{lemma}

Our main result in this section is stated as follows.

\begin{theorem}\label{PoC}
Denote the empirical distribution of the uniform reshuffling stochastic system \eqref{uniform} at time $t$ as \[\rho_{\emp}(t):= \frac{1}{N}\,\sum_{i=1}^N \delta_{X_i(t)},\] and let $q(t)$ be the solution of \eqref{PDE} with initial condition $q(0)$. If \begin{equation}\label{assump}
\mathbb{E}[W_1(\rho_{\mathrm{\emp}}(0),q(0))] \xrightarrow[]{} 0 ~\text{as}~ N\to \infty,
\end{equation}
then we have that\[\mathbb{E}[W_1(\rho_{\mathrm{\emp}}(t),q(t))] \xrightarrow[]{} 0 ~\text{as}~ N\to \infty,\] holding for all $0\leq t \leq T$ with any prefixed $T >0$.
\end{theorem}

\subsection{Switching supremum and expectation}
\label{subsec:5.2}
We will also make use of the following result, which allows us to interchange the operation of supremum and of expectation.
\begin{lemma}\label{exchangesup}
Consider a random Radon measure $Z$ on $\R$ with $\int Z(\dd x)=0$ and with uniformly bounded second moment $\int (1+|x|^2)\,|Z|(\dd x)\leq m_2$ almost surely for some constant $m_2$. Then there exists $\theta>0$ such that
\[\mathbb{E}\,\left[ \sup_{\|\nabla\varphi\|_{\infty}\leq 1} \int \varphi\,\dd Z\right]\leq C\,m_2\,\left(
\sup_{\|\nabla\varphi\|_{\infty}\leq 1} \mathbb{E}\,\left[ \int \varphi\,\dd Z\right]^2\right)^\theta.\]
\end{lemma}

\Proof This is essentially an interpolation argument. First of all, we can always assume that $\varphi(0)=0$ by subtracting a constant.
Introduce a classical convolution kernel $K_\eps$. We have that $\|K_\eps\star \varphi-\varphi\|_{L^\infty}\leq C\,\eps$ which implies that
\[\int \varphi\,Z(dx)\leq \int K_\eps\star\varphi\,Z(\dd x)+C\,\eps.\]
Then we reduce ourselves to a compact support: since $\|\nabla\varphi\|_{\infty}\leq 1$ then $|\varphi(x)|\leq |x|$ and
\[\begin{split}
\int K_\eps\star\varphi\,Z(\dd x)&\leq \int_{|x|\leq R} K_\eps\star\varphi\,Z(\dd x)+2\,\int_{|x|\geq R} |x|\,|Z|(\dd x)\\
&\leq \int_{|x|\leq R} K_\eps\star\varphi\,Z(\dd x)+2\,\frac{m_2}{R}.
\end{split}\]
On $[-R, R]$, we have on the other hand that $\|K_\eps\star\varphi\|_{H^2}\leq \frac{C}{\eps}\,\|\varphi\|_{W^{1,\infty}}\leq C\,\frac{R}{\eps}$. Hence
\[\sup_{\|\nabla\varphi\|_{\infty}\leq 1} \int \varphi\,Z(\dd x)\leq C\,\frac{R}{\eps}\,\sup_{\|\varphi\|_{H^2}\leq 1}\int_{|x|\leq R} \varphi\,Z(\dd x)+C\,m_2\,\left(\eps+\frac{1}{R}\right).\]
Of course \[\sup_{\|\varphi\|_{H^2}\leq 1}\int_{|x|\leq R} \varphi\,Z(\dd x)=\|Z\|_{H^{-2}([-R,R])},\]
and by using Fourier series
\[\|Z\|_{H^{-2}([-R,\ R])}^2=\sum_k \frac{R^2}{1+k^4}\,\left(\int_{-R}^R \expo^{-i\,k\,\pi\,x/R}\,\dd Z\right)^2.\]
Hence by Cauchy-Schwartz,
\[\begin{split}
&\mathbb{E}\left[\sup_{\|\nabla\varphi\|_{\infty}\leq 1} \int \varphi\,Z(\dd x)\right]\leq C\,m_2\,\left(\eps+\frac{1}{R}\right)\\
&\qquad+C\,\frac{R^2}{\eps}\,\left(\sum_k \frac{1}{1+k^4}\,\mathbb{E}\,\left[\left(\int_{-R}^R \expo^{-i\,k\,\pi\,x/R}\,\dd Z\right)^2\right] \right)^{1/2}.
\end{split}\]
Finally we have that
\[\|\nabla \expo^{-i\,k\,\pi\,x/R}\|_{\infty}\leq C\,k,\]
so that
\[\mathbb{E}\,\left[\left(\int_{-R}^R \expo^{-i\,k\,\pi\,x/R}\,\dd Z\right)^2\right]\leq C\,k^2\,\sup_{\|\nabla\varphi\|_{\infty}\leq 1}\mathbb{E}\,\left[\left(\int \varphi\,\dd Z\right)^2\right].\]
This allows us to conclude that
\[\begin{split}
&\mathbb{E}\left[\sup_{\|\nabla\varphi\|_{\infty}\leq 1} \int \varphi\,Z(\dd x)\right]\leq C\,m_2\,\left(\eps+\frac{1}{R}\right)\\
&\qquad+C\,\frac{R^2}{\eps}\,\left(\sum_k \frac{k^2}{1+k^4}\,\sup_{\|\nabla\varphi\|_{\infty}\leq 1}\mathbb{E}\,\left[\left(\int \varphi\,\dd Z\right)^2\right] \right)^{1/2},
\end{split}\]
or \[\begin{split}
&\mathbb{E}\left[\sup_{\|\nabla\varphi\|_{\infty}\leq 1} \int \varphi\,Z(\dd x)\right]\leq C\,m_2\,\left(\eps+\frac{1}{R}\right)\\
&\qquad+C\,\frac{R^2}{\eps}\,\left(\sup_{\|\nabla\varphi\|_{\infty}\leq 1}\mathbb{E}\,\left[\left(\int \varphi\,\dd Z\right)^2\right] \right)^{1/2},
\end{split}\]
which finishes the proof by optimizing in $R$ and $\eps$.\qed

\subsection{Proof of propagation of chaos}
\label{subsec:5.3}
\noindent The proof of Theorem \ref{PoC} occupies the rest of the section. \\
\Proof We recall that the map $Q_+[\cdot] \colon \mathcal{P}(\mathbb{R}_+) \to \mathcal{P}(\mathbb{R}_+)$ is defined via
\[Q_+[q](x) = \int_0^\infty\int_0^\infty \frac{\mathbbm{1}_{[0,k+\ell]}(x)}{k+\ell}q(k)q(\ell)\,\dd k\, \dd \ell,\] and that a classical solution $q(t,x)$ of
\begin{equation}\label{APDE}
q(t,x) = q(0,x) + \int_0^t G[q](s,x)\,\dd s
\end{equation}
exists for $0\leq t < \infty$, where $G = Q_+ - \mathrm{Id}$ and $q(0,x)$ is an continuous probability density function with mean $m_1$ whose support is contained in $\mathbb{R}_+$. The map $Q_+$ is Lipschitz continuous in the sense that
\begin{equation}\label{LipG}
W_1(Q_+[f],Q_+[g]) \leq W_1(f,g)
\end{equation}
for any $f,g \in \mathcal{P}(\mathbb{R}_+)$. Indeed, we have \[W_1(Q_+[f],Q_+[g]) = \sup\limits_{\|\nabla \varphi\|_\infty \leq 1} \mathbb{E}\left[\varphi(U(X_1+Y_1)) - \varphi(U(X_2+Y_2)\right],\] where $X_1, Y_1$ are i.i.d with law $f$, $X_2, Y_2$ are i.i.d with law $g$, and $U \sim \mathrm{Uniform}[0,1]$ is independent of $X_i$ and $Y_i$ for $i=1,2$. By Lipschitz continuity of the test function $\varphi$, we obtain
\[W_1(Q_+[f],Q_+[g]) \leq \mathbb{E}\left[2\,U|X_1-X_2|\right] = \mathbb{E}[|X_1-X_2|].\] We now recall an alternative formulation of $W_1(f,g)$, given by \[W_1(f,g) = \inf\left\{\mathbb E[|X-Y|];~\mathrm{Law}(X)=f,~\mathrm{Law}(Y)=g\right\},\] so in particular, we may take a coupling of $X_1$ and $X_2$ so that $W_1(f,g) = \mathbb{E}[|X_1-X_2|]$. Assembling these pieces together, we arrive at \eqref{LipG}.

We are going to prove a more precise control than \eqref{LipG}, by working directly on $Q_+[f]$. Consider now two random probability measures $f$ and $g$ with bounded second moment and a deterministic test function $\varphi$. We have that
\begin{align*}
\int \varphi(x)\,(Q_+[f]-Q_+[g])\, \dd x &=\int \frac{\mathbbm{1}_{x\leq k+\ell}}{k+\ell} \varphi(x)\,(f(\dd k)-g(\dd k))\,(f(\dd \ell)+g(\dd \ell))\,\dd x\\
&=\int (f(\dd \ell)+g(\dd \ell))\,\int \Phi_{\ell}(k)\,(f(\dd k)-g(\dd k)),
\end{align*}
where we denote
\[\Phi_{\ell}(k)=\frac{1}{k+\ell}\,\int_0^{k+\ell} \varphi(x)\,\dd x.\]
Since $\int Q_+[f]=\int Q_+[g]$, we can always assume without loss of generality that $\varphi(0)=0$, whence $|\varphi(x)|\leq \|\nabla\varphi\|_{\infty}\,|x|\leq |x|$.
Now we observe that $\Phi_{\ell}$ is deterministic with
\begin{equation}
|\partial_k \Phi_{\ell}(k)|\leq \frac{|\varphi(k+\ell)|}{k+\ell}+\frac{1}{(k+\ell)^2}\,\int_0^{k+\ell} |\varphi(x)|\,\dd x\leq 1+\frac{1}{(k+\ell)^2}\,\int_0^{k+\ell} x\,\dd x\leq \frac{3}{2}.\label{nablaPhil}
\end{equation}
By \eqref{nablaPhil} and recalling that again $\Phi_{\ell}$ is deterministic and obtained from $\varphi$, we obtain:
\[\mathbb{E} \left[\int \Phi_{\ell}(k)\,(f(\dd k)-g(\dd k))\right] \leq \frac{3}{2}\,\mathbb{E} \left[\sup_{\|\nabla\varphi\|_{\infty}\leq 1}\int \varphi(x)\,(f(\dd x)-g(\dd x))\right].
\]
Therefore we conclude that
\begin{equation}\label{newLipG}
\mathbb{E} \left[\sup_{\|\nabla\varphi\|_{\infty}\leq 1} \int \varphi(x)\,(Q_+[f]-Q_+[g])\right]
\leq 3\,\mathbb{E} \left[\sup_{\|\nabla\varphi\|_{\infty}\leq 1} \int \varphi(x)\,(f(\dd x)-g(\dd x))\right].
\end{equation}
We now observe that the empirical measure is a compound jump process: Define $N_t$ a homogeneous Poisson process with constant intensity $\lambda=(N-1)/2$. Given $\tau_1,\ldots,\tau_k$ the times when $N_t$ jumps, we take the $Y_{\tau_k}$ independent: At each $\tau_k$, with uniform probability $\frac{2}{N\,(N-1)}$ we choose a pair $i< j$ and take
\[\begin{split}
Y_{\tau_k}=&\frac{1}{N}\,\Big(\delta(x-U_{k}\,(X_i(\tau_k-)+X_j(\tau_k-))+\delta(x-(1-U_{k})\,(X_i(\tau_k-)+X_j(\tau_k-))\\ &-\delta(x-X_i(\tau_k-))-\delta(x-X_j(\tau_k-))\Big),
\end{split}
\]
where the $U_k$ are i.i.d. in $[0,1]$.

We immediately note that
\begin{equation}
\begin{split}
&\lambda\,\mathbb{E}[Y_t]=\frac{1}{N^2}\,\sum_{i<j} \mathbb{E}\Big[\delta(x-U\,(X_i(t-)+X_j(t-))\\&\ +\delta(x-(1-U)\,(X_i(t-)+X_j(t-)) -\delta(x-X_i(t-))-\delta(x-X_j(t-))\Big],
\end{split}\label{expectYt}
\end{equation}
where $U$ is uniformly distributed in $[0,1]$ and independent of all $X_i(t-)$.

We also remark that by the standard control of moments, we immediately have that
\begin{equation}
\int x^2\,\rho_{\emp}(t,\dd x)\leq m_2=2+\int x^2\,\rho_{\emp}(0,\dd x).\label{secondmomentempiric}
  \end{equation}
We now show that the empirical measure of the stochastic system satisfies an approximate version of \eqref{APDE}. Fix a deterministic test function $\varphi$ with $\|\nabla \varphi\|_\infty \leq 1$, and consider the time evolution of $\langle \rho_{\emp},\varphi \rangle$ where for some probability measure $\nu$, we denote by the duality bracket $\langle \nu,\varphi \rangle =\int \varphi\,\dd \nu$. We emphasize here that $\varphi$ can also be random and will indeed be chosen according to $\rho_{\emp}$ to estimate Wasserstein distances involving $\rho_{\emp}$. Then
\[\dd \mathbb{E}[\langle \rho_{\emp},\varphi \rangle]=\dd \mathbb{E}\left[\langle Y_t\,\dd N_t,\varphi \rangle \right] =\lambda\,\langle\mathbb{E}[Y_t],\varphi \rangle\,\dd t.\]
Hence by \eqref{expectYt},
\[
\begin{aligned}
\dd \mathbb{E}[\langle \rho_{\emp},\varphi \rangle] &= \frac{1}{N^2}\,\sum_{i < j} \mathbb{E}\left[\varphi\big(U(X_i+X_j)\big) + \varphi\big((1-U)(X_i+X_j)\big) - \varphi(X_i) - \varphi(X_j)\right] \dd t \\
&=\frac{1}{N^2}\,\sum_{i,j=1\ldots N,i\neq j} \mathbb{E}\left[\varphi\big(U(X_i+X_j)\big)-\varphi(X_i)\right] \dd t\\
&=\frac{1}{N^2}\,\sum_{i,j=1}^N \mathbb{E}\left[\varphi\big(U(X_i+X_j)\big)-\varphi(X_i)\right] \dd t + R\,\dd t,
\end{aligned}
\]
where all $X_i,\,X_j$ are taken at time $t-$ and where $R=-\frac{1}{N^2}\,\sum_i \mathbb{E}\left[\varphi\big(2\,U\,X_i\big)-\varphi(X_i)\right]$. Hence $|R| \leq \mathcal{O}\left(\frac{1}{N}\right)$ uniformly over $\varphi$ and $t\geq 0$.
On the other hand, we may calculate
\[\langle Q_+[\rho_{\emp}], \varphi \rangle=\frac{1}{N^2}\,\sum_{i,j}\int \varphi(x)\,\frac{\mathbbm{1}_{x\leq X_i+X_j}}{X_i+X_j}\,\dd x=\frac{1}{N^2}\,\sum_{i,j}\int_0^1 \varphi(u\,(X_i+X_j))\,\dd u,\]
by the change of variables $x=u\,(X_i+X_j)$. Therefore
\begin{equation}
\dd \mathbb{E}[\langle \rho_{\emp},\varphi \rangle]=\mathbb{E}\left[\langle G[\rho_{\emp}], \varphi \rangle\right]\,\dd t+R\,\dd t.\label{weakformulation}
\end{equation}
By Dynkin's formula, the compensated process
\begin{equation}\label{mart}
M_\varphi(t):= \langle \rho_{\emp}(t),\varphi \rangle - \langle \rho_{\emp}(0),\varphi \rangle - \int_0^t \left(\mathbb{E}[\langle G[\rho_{\emp}(s)], \varphi \rangle] + R(s)\right)\, \dd s
\end{equation}
is a martingale. Furthermore, comparing with \eqref{APDE}, we easily obtain that
\[\begin{aligned}
\langle \rho_{\emp}(t) - q(t),\varphi \rangle &= M_\varphi(t) + \langle \rho_{\emp}(0) - q(0),\varphi \rangle \\
&\quad + \mathbb{E}\,\int_0^t \langle G[\rho_{\emp}(s)] - G[q(s)], \varphi \rangle \,\dd s + \mathcal{O}\left(\frac{t}{N}\right).
\end{aligned}\]
Taking the supremum over $\varphi$, we therefore have that
\[\begin{aligned}
&\mathbb{E}\sup_{\|\nabla\varphi\|_{\infty}\leq 1} \langle \rho_{\emp}(t) - q(t),\varphi \rangle \leq \mathbb{E}\,\sup_{\|\nabla\varphi\|_{\infty}\leq 1} (|M_\varphi(t)| + \langle \rho_{\emp}(0) - q(0),\varphi \rangle) \\
&\quad +\int_0^t \mathbb{E}\,\sup_{\|\nabla\varphi\|_{\infty}\leq 1} \langle G[\rho_{\emp}(s)] - G[q(s)], \varphi \rangle \,\dd s + \mathcal{O}\left(\frac{t}{N}\right).
\end{aligned}\]
By the definition of the $W_1$ distance, we deduce from \eqref{newLipG}  that
\[\begin{split}
&\mathbb{E}\,W_1(\rho_{\emp}(t),q(t)) \leq \eta(t) + C\,\int_0^t \mathbb{E}\, W_1(\rho_{\emp}(t),q(t))\, \dd s + \frac{C\,t}{N},
\end{split}\]
in which we have set
\begin{equation}\label{eta}
\eta(t) :=  \mathbb{E}\,\sup\limits_{\|\nabla \varphi\|_\infty \leq 1}\,|M_\varphi(t)| + \mathbb{E}\,W_1(\rho_{\emp}(0),q(0)).
\end{equation}
Thus, Gronwall's inequality gives rise to
\begin{equation}\label{almostfinish}
\mathbb{E}\,W_1(\rho_{\emp}(t),q(t))\leq \left(\sup\limits_{t\in [0,T]} \eta(t) + \frac{C\,T}{N}\right)\expo^{C\,T}.
\end{equation}
In order to establish propagation of chaos for $t \leq T$, it therefore suffices to show that
\begin{equation}\label{convprob}
\sup\limits_{t\in [0,T]} \eta(t) \xrightarrow[]{N \to \infty} 0.
\end{equation}
To prove \eqref{convprob}, we treat each term appearing in the definition of $\eta(t)$ separately. The second term in \eqref{eta} approaches to 0 as $N \to \infty$ by our assumption.

To handle the first term, let us write $Z(t) = \langle \rho_{\emp}(t),\varphi \rangle$ and $M(t) = M_\varphi(t)$ for notation simplicity. Of course $Z(t)$ is a compound jump process itself and by combining \eqref{weakformulation} and \eqref{mart}
\[M_\varphi(t)=Z(t)-Z(0)-\int_0^t \tilde Y(s)\,ds,\quad \tilde Y(t)=\langle G[\rho_{\emp}(t)],\varphi\rangle +R.\]
We may hence use It\^o's lemma as stated in Lemma \ref{Ito}, which yields
\[
\dd \mathbb{E}[M^2(t)] = \sum_{i<j} \mathbb{E}\left[M^2_{ij}(t) - M^2(t)\right]\,\frac{\dd t}{N} - \mathbb{E}\left[2\,M(t)\,\langle G[\rho_{\emp}(t)],\varphi \rangle\right]\,\dd t + \mathcal{O}\left(\frac{1}{N}\right)\dd t,
\]
where $M_{ij} = M + Y_{ij}$ and we define
\[
Y_{ij} := \big\langle \frac{1}{N}(\delta_{U_k(X_i+X_j)} + \delta_{(1-U_k)\,(X_i+X_j)} - \delta_{X_i} - \delta_{X_j}), \varphi \big\rangle.
\]
Therefore, we have
\begin{equation*}
\begin{aligned}
  \dd \mathbb{E}[M^2(t)] =& \sum_{i<j} \mathbb{E}\left[2\,M(t)\,Y_{ij} + Y^2_{ij} \right]\,\frac{\dd t}{N} - \mathbb{E}\left[2\,M(t)\,\langle G[\rho_{\emp}(t)],\varphi \rangle\right]\,\dd t\\
  &+ \mathcal{O}\left(\frac{1}{N}\right)\dd t. \\
\end{aligned}
\end{equation*}
By our previous calculations
\[
\begin{split}
  &\frac{1}{N}\,\sum_{i<j} \mathbb{E}[M(t)\,Y_{ij}]\\
  &\quad=\frac{1}{N^2}\,\sum_{i<j} \mathbb{E}[M(t)\, (\varphi(U\,(X_i+X_j) +\varphi((1-U)\,(X_i+X_j) -\varphi(X_i) -\varphi(X_j))] \\
  &\quad=\frac{1}{N^2}\,\sum_{i\neq j} \mathbb{E}[M(t)\,(\varphi(U\,(X_i+X_j))-\varphi(X_i))]\\
  &\quad=\frac{1}{N^2}\,\sum_{i,j} \mathbb{E}[M(t)\,(\varphi(U\,(X_i+X_j))-\varphi(X_i))]+O\left(\frac1N\right),
\end{split}
\]
as $U$ is random variable independent of $M(t)$ and $\rho_{emp}(t)$.

Therefore
\[
\frac{1}{N}\,\sum_{i<j} \mathbb{E}[M(t)\,Y_{ij}]=\mathbb{E}\left[M(t)\,\langle G[\rho_{\emp}(t)],\varphi \rangle\right]+O\left(\frac1N\right),
\]
and consequently
\[\dd \mathbb{E}[M^2(t)] =\sum_{i<j} \mathbb{E}\left[Y^2_{ij} \right]\,\frac{\dd t}{N} +\mathcal{O}\left(\frac{1}{N}\right)\dd t\leq \frac{C}{N}\dd t, \\\]
for a constant $C$ that depends only on $\|\nabla\varphi\|_\infty$.
This lets us deduce that
\[\sup_{\|\nabla\varphi\|_\infty\leq 1}\,\mathbb{E}\,\left[M_\varphi(t)\right]\leq \frac{C\,t}{N}.\]
Recalling the definition of $M_\varphi(t)$, we have that
\[M_\varphi(t)=\int \varphi(x)\,\mu(t,\dd x)\]
for some random Radon measure $\mu$ with uniformly bounded second moment. Furthermore $\int \mu(t,\dd x)=0$ since $\int \rho_{\emp}(t,\dd x)=1=\int \rho_{\emp}(0,\dd x)$ and $\int G[\rho_{\emp}(t)]\,\dd x =0$.

\noindent We may hence apply Lemma \ref{exchangesup} to obtain that
\[\mathbb{E}\,\left[\sup_{\|\nabla\varphi\|_\infty\leq 1}\, M_\varphi(t)\right]\leq C\,\frac{t^\theta}{N^\theta},\]
which allows to conclude that $\sup\limits_{t\in [0,T]} \eta(t) \xrightarrow[]{N\to\infty} 0$. \qed

\begin{remark}
One can readily check that \[\|Q_+[f] - Q_+[g]\|_{L^1(\mathbb{R}_+)} \leq 2\,\|f - g\|_{L^1(\mathbb{R}_+)}\] for all probability densities $f,g$ whose support are contained in $\mathbb{R}_+$, but as we are working on $\mathcal{P}(\mathbb{R}_+)$, we can not use any strong distances. Hence, equipping $\mathcal{P}(\mathbb{R}_+)$ with an appropriate distance so that the operator $Q_+$ has enjoys a Lipschitz continuity with respect to the chosen distance is an indispensable step to make the argument above work.
\end{remark}

\section{Appendix}

\subsection{Proof of Theorem \ref{thm3}}
\Proof The whole strategy is of course to find some $\delta$ such that if
\begin{equation}
  \int q(t,x)\,\log \frac{q(t,x)}{q_\infty(x)}\,\dd x\leq \delta,\label{smallentropy}
\end{equation}
then we have for the $\varepsilon$ of Theorem~\ref{localconvergence}
\begin{equation}
\int \frac{|q(t,x)-q_\infty(x)|^2}{q_\infty(x)}\,dx\leq\varepsilon.\label{smallL2}
  \end{equation}
We start with using Lemma \ref{label1} for $C=2$ and note that
\begin{equation}
\begin{split}
  &\frac{1}{4}\,\int_{q_\infty/2\leq q\leq 2\,q_\infty} \frac{|q(t,x)-q_\infty(x)|^2}{q_\infty(x)}\,\dd x+\frac{1}{8}\,\int_{q\leq q_\infty/2} q_\infty(x)\,\dd x+\frac{\log 2}{4}\int_{q\geq 2\,q_\infty} q(t,x)\,\dd x\\
  &\qquad\leq \int q(t,x)\,\log \frac{q(t,x)}{q_\infty(x)}\,\dd x.
\end{split}\label{lowerboundl2entropy}
\end{equation}
Observe that if $q\leq q_\infty/2$ then
\[
\frac{|q(t,x)-q_\infty(x)|^2}{q_\infty(x)}\leq q_\infty(x),
\]
so the first two terms already provides the straightforward bound
\begin{equation}\label{l2entropyfirst}
\int_{q\leq 2\,q_\infty} \frac{|q(t,x)-q_\infty(x)|^2}{q_\infty(x)}\,\dd x\leq 8\,\int q(t,x)\,\log \frac{q(t,x)}{q_\infty(x)}\,\dd x.
\end{equation}
Now if $q\geq q_\infty$ then
\[\frac{|q(t,x)-q_\infty(x)|^2}{q_\infty(x)}\leq \frac{(q(t,x))^2}{q_\infty(x)}.\]
Therefore for any $p > 1$,
\begin{equation*}
\begin{split}
  &\int_{q\geq 2\,q_\infty} \frac{|q(t,x)-q_\infty(x)|^2}{q_\infty(x)}\,\dd x\leq \int_{q\geq 2\,q_\infty} \frac{|q(t,x)|^2}{q_\infty(x)}\,\dd x\\
  &\qquad\leq \left(\int_{q\geq 2\,q_\infty} q(t,x)\,\dd x\right)^{1-1/p}\,\left(\int_{q\geq 2\,q_\infty} \frac{|q(t,x)|^{p+1}}{\left(q_\infty(x)\right)^p}\,\dd x\right)^{1/p}.
\end{split}
\end{equation*}
We now use Corollary~\ref{exponentialpointwise} to find that
\[
\begin{split}
  &\int_{q\geq 2\,q_\infty} \frac{|q(t,x)|^{p+1}}{\left(q_\infty(x)\right)^p}\,\dd x\leq C_p\, \int \left(\expo^{(-(p+1)\,\lambda_0+p)\,x}+\expo^{p\,x}\,(q(0,x))^{p+1}\right)\,\dd x\\
  &\qquad\leq C_p'\,\int \expo^{(-(p+1)\,\lambda + p)\,x}\,\dd x,
\end{split}
\]
in which $\lambda \in (\frac 12, \lambda_0)$. Now we take $p$ close enough to $1$ such that $p-(p+1)\,\lambda < 0$ which is always possible if $\lambda_0>\frac 12$. For this choice of $p$, we hence obtain that
\[
\int_{q\geq 2\,q_\infty} \frac{|q(t,x)-q_\infty(x)|^2}{q_\infty(x)}\,\dd x\leq C_p\,\left(\int_{q\geq 2\,q_\infty} q(x)\,\dd x\right)^{1-1/p}.
\]
Going back to \eqref{lowerboundl2entropy}, we can conclude that
\[\int_{q\geq 2\,q_\infty} \frac{|q(t,x)-q_\infty(x)|^2}{q_\infty(x)}\,\dd x\leq C_p\,\left(\int q(t,x)\,\log\frac{q(t,x)}{q_\infty(x)}\,\dd x\right)^{1-1/p},\]
and combining this with \eqref{l2entropyfirst}, we deduce that for some $C$ and $\theta \in (0,1)$
\[
\int \frac{|q(t,x)-q_\infty(x)|^2}{q_\infty(x)}\,\dd x\leq C\,\left(\int q(t,x)\,\log\frac{q(t,x)}{q_\infty(x)}\,\dd x\right)^{\theta}\leq C\,\delta^{\theta}.
\]
It is enough to choose $\delta$ being small enough to conclude the proof. \qed

\bibliographystyle{plain}
\bibliography{uniform_reshuffling}

\end{document}